\documentclass[11pt, reqno]{amsart}
\usepackage{amsmath,amsfonts,amssymb,amsthm}
\usepackage{epsfig}

\usepackage{amstext}
\usepackage{xspace}
\usepackage{amsfonts}
\usepackage{amsmath}
\usepackage{amssymb}
\usepackage{amstext}
\usepackage{amsthm}       %proof environment :)
\usepackage{xspace}

\def\supess{\mathop{\rm ess\: sup }}
\newcommand{\Sf}{\mathbb S}
\newcommand{\C}{\mathbb C}
\newcommand{\EE}{\mathbb E}
\newcommand{\R}{\mathbb R}
\newcommand{\Length}{\mathbb L}
\newcommand{\N}{\mathbb N}
\newcommand{\Q}{\mathbb Q}
\newcommand{\T}{\mathbb T}
\newcommand{\Z}{\mathbb Z}
\newcommand{\restr}{\hbox{\LARGE $\llcorner$}}
\newcommand{\LL}{\mbox{\restr}}
\newcommand{\comp}{\mbox{\scriptsize  $\circ$}}
\newcommand{\eps}{\varepsilon}
\newcommand{\tto}{\longrightarrow}
\newcommand{\ucv}{\rightrightarrows}

\newcommand{\A}{\mathcal{A}}
\newcommand{\hh}{\mathcal{H}}
\newcommand{\ddiv}{\textrm{\rm div}}
\newcommand{\e}{\textrm{\rm e}}
\newcommand{\Eq}{\mathcal{E}}
\newcommand{\F}{\mathcal{F}}
\newcommand{\K}{\mathcal{K}}
\newcommand{\M}{\mathcal{M}}
\renewcommand{\S}{\mathcal{S}}
\newcommand{\SLip}{L^\infty(\Omega;\D{Lip}(\R^N))}
\newcommand{\leb}{\mathcal{L}}
\newcommand{\lebbar}{\bar{\mathcal{L}}}
\newcommand{\dist}{\mathcal{D}}
\newcommand{\massa}{\mbox{\boldmath{$M$}}}
\newcommand{\PP}{\mathbb{P}}
\newcommand{\tagliato}{$\kern-5.5 mm -$}
\newcommand{\tagliat}{$\kern-6.5 mm -$}
\newcommand{\taglio}{\mbox{\tagliato}}
\newcommand{\tagli}{\mbox{\tagliat}}
\newcommand{\cchi}{\mbox{\large $\chi$}}
\newcommand{\abra}[1]{(\ref{#1})}
\newcommand{\B}[1]{\mbox{\boldmath $#1$}}
\newcommand{\D}[1]{\mbox{\rm #1}}
\newcommand{\dd}{\D{d}}
\newcommand{\cv}[1]
{\, \displaystyle{\mathop{\longrightarrow}\limits_{#1}}\, }

\newtheorem{teorema}{Theorem}[section]
\newtheorem{prop}[teorema]{Proposition}
\newtheorem{lemma}[teorema]{Lemma}
\newtheorem{definition}[teorema]{Definition}

\newtheorem{guess}[teorema]{Remark}
\newtheorem{example}[teorema]{Example}

\newenvironment{esempio}{\begin{example} \begin{rm}}{\end{rm} \end{example}}
\newenvironment{dimo}{{\bf\noindent Proof.}}{\qed}
\newenvironment{oss}{\begin{guess} \begin{rm}}{\end{rm} \end{guess}}
\newenvironment{definizione}{\begin{definition} \begin{rm}}{\end{rm}
\end{definition}}
\begin{document}

\title{A metric analysis\\ of critical  Hamilton--Jacobi equations\\
in the stationary ergodic setting}
\author{Andrea Davini \and Antonio Siconolfi}
\address{Dip. di Matematica, Universit\`a di Roma ``La Sapienza",
P.le Aldo Moro 2, 00185 Roma, Italy}
\email{davini@mat.uniroma1.it, siconolf@mat.uniroma1.it}
% \address{Dip. di Matematica, Universit\`a di Roma ``La Sapienza",
% P.le Aldo Moro 2, 00185 Roma, Italy}
% \email{siconolf@mat.uniroma1.it}
%\date{\today}
%\subjclass{41A99, 49J45, 58B20}

\begin{abstract}  We adapt the metric approach to the study of  stationary ergodic
Hamilton--Jacobi equations, for which a notion of admissible
random (sub)solution is defined. For any level of the Hamiltonian
greater than or equal to a distinguished critical value, we define
an intrinsic random semidistance and prove that an asymptotic norm
does exist. Taking as source region a suitable class of closed
random sets, we show  that the Lax formula provides admissible
subsolutions. This enables us to relate the degeneracies of the
critical stable norm to the existence/nonexistence of exact or
approximate critical admissible solutions.

\end{abstract}
\maketitle

\begin{section}{Introduction} \label{over}

The main purpose of the paper is to adapt the so--called metric
method, which has revealed to be a powerful tool for the analysis
of  critical Hamilton--Jacobi equations posed on compact spaces,
see \cite{CS2, DavSic05, FaSic03}, to the stationary ergodic
setting. Loosely speaking,  the ergodicity can be viewed as  a
weaker form of compactness,  mostly thanks to some powerful
asymptotic results, like Birkhoff and Kingman subadditive
Theorems, that we repeatedly employ in our research.

We consider a probability space $\Omega$, on which the action of
$\R^N$  gives rise to an $N$--dimensional ergodic dynamical
system, and a random continuous  Hamiltonian $H(x,p,\omega)$,
which is stationary with respect to such dynamics, and, in
addition, convex and coercive in the momentum variable.\par

We look for {\em admissible subsolutions} of the corresponding
stochastic Hamilton--Jacobi equations at  different levels of the
Hamiltonians. By this we mean Lipschitz random functions which are
almost surely subsolutions either in the {\em viscosity} sense,
or, equivalently, {\em almost everywhere}, while
 the term {\em admissible} refers to the fact that they are stationary or,
 in a weaker form, that they
possess stationary increments and gradient with vanishing mean.
Exploiting  ergodicity and Birkhoff Theorem, this last property
turns out to be equivalent to the almost sure sublinearity  at
infinity.\par

Actually, we prove that the infima of the values  for which the
corresponding equations admit a subsolution of the two types
coincide. This quantity is called the {\em stationary critical
value} of $H$ and will be denoted by $c$. The difference is that,
due to lack of stability, the existence of a stationary
{subsolution} to $H(x,D v,\omega)=c$ can fail already in the
one--dimensional setting, see \cite{DavSic08}, while an
Ascoli--type theorem, see Theorem \ref{stabile}, adjusted to the
random environment, guarantees  to find subsolutions of the latter
class to   the critical equation. For this reason we will use,
from now on, the word admissible in the weak sense.\par

The relevance of the critical value is in the fact that it is the
unique level of $H$ for which the corresponding Hamilton--Jacobi
equation can have admissible (viscosity) solutions or approximate
solutions, see Section \ref{sez HJ} for definitions.  These
objects can be used as exact or approximate correctors in related
homogenization procedures implementing the perturbed test function
method \cite{Ev89, Ev92}.\par

Existence and nonexistence issues for exact and approximate
solutions are relevant open problems in the field. So far, the
setup has been completely clarified only in the one--dimensional
case \cite{DavSic08}, where we proved the existence of approximate
or exact correctors, depending on whether $0$ belongs or not to
the interior of the flat part of the effective Hamiltonian.
%
%, where it has been proved the existence of exact correctors
%outside the interior of the flat part of the effective Hamiltonian
%and that of approximate correctors inside this region.

What is disappointing, at first sight,  about the   metric
approach in this context, is that it is, in the starting point,
purely deterministic, with $\omega$ playing just the role of a
parameter.  We in fact define for every $\omega \in \Omega$ and $a
\in \R$ an intrinsic  semidistance $S_a$ starting from  the
support function of the $a$--sublevel of $H(x,\cdot, \omega)$.
\par

 It is
well known that such a distance is finite if and only if
$H(x,Du,\omega)=a$ admits (deterministic) subsolutions. A new
critical value, depending on $\omega$, say $c_f(\omega)$, then
comes at the surface, corresponding to the minimum $a$ for which
the equation admits subsolutions. Because of the measurability
properties of the Hamiltonian the map $\omega \mapsto c_f(\omega)$
is a random variable, which is, in addition, almost surely
constant by the stationarity of $H$ and the ergodicity assumption.
Such a constant will be denoted by $c_f$ and called {\em free
critical value} to distinguish it from $c$. For the same reasons,
$S_a$ is, for $a \geq c_f$,  a stochastic semidistance, namely a
random variable taking value in the family of semidistances
endowed with the local uniform convergence. It is apparent from
its very definition that $c_f \leq c$, and strict inequality is
possible.\par

  From what previously outlined, it is clear that the
intrinsic distances $S_a$ cannot be useful {\em per se} to our
analysis, in particular the critical level cannot be detected
through the appearance of some degeneracies  of the corresponding
intrinsic  random distance,  like in the compact setting. As a
matter of fact, such kind of phenomena do not take place, in
general, even when $a < c$. Some other steps should  therefore be
accomplished.
\par

We basically follow two ways: first, we perform an asymptotic
analysis  of  intrinsic distances showing that corresponding
(deterministic) {\em stable norms}, say $\phi_a$, do exist for any
$a \geq c_f$  and enjoy some relevant properties; {secondly}, we
generalize Lax--type formulae to the stochastic environment
providing a class of admissible subsolutions. Through the
interplay of these lines of investigation, we establish  in the
end some of our main results.\par

To show that there is   a stable norm, even in the deterministic
case,  some kind of subadditive principle is needed, see for
instance \cite{Bu}. Here we use  Kingman's Subadditive Ergodic
Theorem and mimic the proof given in \cite{ReTa00, Souga99} for
the existence of an {\em effective} (homogenized) {\em
Hamiltonian}. The stable norms $\phi_a$ are of Minkowski type,
possibly degenerate; we actually prove that the critical value is
the infimum of the $a$ for which $\phi_a$ is nondegenerate. \par

Such norms, being {convex} and positively homogeneous, are the
support functions of some compact convex  sets which  can be
interpreted  as the associate dual unit balls.  We show that they
coincide with the corresponding sublevels of the effective
Hamiltonian, denoted by $\overline H$. In this manner we provide a
new simpler proof of a result already established in
\cite{LiSou03} through PDE techniques, namely that the effective
Hamiltonian coincides with the function  associating with any $P
\in \R^N$ the stationary critical value of the Hamiltonian
$H(x,P+p,\omega)$, see  Theorem \ref{effettivo}. Moreover we show that
the free critical value is the minimum of $\overline H$;  an
analogous result have been obtained in \cite{FatMad} for
Hamiltonians defined in unbounded spaces and enjoying some form of
symmetry.

Regarding the Lax formula, we recall that  for any fixed $\omega$
a class of fundamental subsolutions to $H(x,Du,\omega)=a$ is built
up by \[\inf\{g(y)+ S_a(y,x,\omega)\,:\, y \in C\}, \] where $C$
is a closed subset and $g$  a function defined on it which is
$1$--Lipschitz continuous with respect to $S_a$. These kinds of
functions are, in addition, solutions outside $C$. To get through
this pattern {\em admissible} subsolutions, appropriate conditions
have to be assumed on the source set, which depends on $\omega$,
as well as on the trace, linking them to the stationary ergodic
structure.
\par

The key idea, already exploited in \cite{DavSic08}, is to borrow
some tools from stochastic geometry (see \cite{Molchanov} for a
comprehensive treatment of this topic), and to take as source
region a stationary closed random set. That is to say a random
variable taking values in the family of closed subsets of $\R^N$
endowed with the Fell topology which, in addition, satisfies a
compatibility property with the ergodic dynamics, see \eqref{def
stationary set}.\par

With this choice the Lax formula gives an
admissible subsolution for any stationary Lipschitz random
function $g$,  provided it takes finite values, see Proposition
\ref{prop Lax in S}. This latter condition is always fullfilled when
$g$ is itself an admissible subsolution, see Proposition \ref{prop Lax in S0}.
In this instance,
the more delicate item to prove is that the function so obtained
is sublinear at infinity, and for this scope it is essential the
asymptotic formula for random closed stationary sets which, in
turn, relies upon Birkhoff Theorem. \par

We use the information gathered to investigate on the existence
of exact correctors when $c=c_f$. At this level, some degeneracies
of the intrinsic semidistance may appear. The collection
of points around which the latter fails to be equivalent to the Euclidean
one form the classical {\em Aubry set}  $\A_f(\omega)$, which, in this setting,
turns out to be  closed random and stationary.
Thus, when it
is almost surely nonempty,  $\A_f(\omega)$ can be used as source region in the Lax formula
to construct an exact corrector. If, on the contrary, it is almost surely empty
and, in addition, the stable norm $\phi_c$ is nondegenerate, in other terms if
no metric degeneracies take place at finite points or at infinity,
then no correctors can exist, see Theorem \ref{teo correctors}.
Note that all known counterexamples to the existence of correctors
are in this frame.\par

When $c=c_f$ and the latter  agrees with $\sup_x
\inf_p H(x,p,\omega)$ almost surely, we also prove that
approximate correctors
can be constructed, always exploiting Lax formula, by taking as source
region the set of $\delta$--approximate equilibria, see Proposition
\ref{equi}. It should be interesting to prove or disprove that
such a property holds true whenever $c=c_f$\par
%
%
% We show that the classical Aubry set $\A_f(\omega)$, made up by
% points around which the intrinsic distance at the free critical
% value fails to be equivalent to the Euclidean one, is closed
% random stationary. Therefore, if it is almost surely nonempty and
% $c=c_f$, a corrector can be constructed via the Lax formula taking
% $\A_f(\omega)$ as  source region. On the contrary, if it is empty
% and the stable norm $\phi_c$ is nondegenerate, in other terms if
% no metric degeneracies take place at finite points or at infinity,
% then no correctors can exist, see Theorem \ref{teo correctors}.
% Note that all known counterexamples to the existence of correctors
% are in this frame.\par
% We also prove that  if $c=c_f$ and the latter  agrees with $\sup_x
% \inf_p H(x,p,\omega)$ almost surely, then approximate correctors
% can be constructed, always exploiting Lax formula, see Proposition
% \ref{equi}. It should be interesting to prove or disprove that
% such a property holds true whenever $c=c_f$.
%\medskip\par

The paper is organized as follows: in Section \ref{pre}  we fix
notations and expose some preliminary material. In Section
\ref{statvar} we introduce two different definitions of
measurability for set--valued variables and the notion of
stationarity; further we describe the main properties of the class
of admissible random functions. Some proofs are postponed to the
Appendix. Section \ref{sez HJ} is focused on stochastic
Hamilton--Jacobi equations and Lax formulae. In Section
\ref{stable} we show the existence of the stable norms associated
with the intrinsic distances and we study their
connection with the effective Hamiltonian. Section \ref{main} is
devoted to statements and proofs of our final  results, and ends  with
an example showing up some major differences with the
one--dimensional setting.

\bigskip

\indent{\textsc{Acknowledgements. $-$}} The first author has been
supported for this research by the European Commission through a
Marie Curie Intra--European Fellowship, Sixth Framework Program
(Contract MEIF-CT-2006-040267). He wishes to thank Albert Fathi
for many interesting discussions and suggestions.\par
%
%\vskip 2cm
\end{section}

\begin{section}{Preliminaries}\label{pre}

We write below a list of symbols used throughout this paper.
\[
\begin{array}{ll}
N & \hbox{an integer number}\\
B_R(x_0) &\hbox{the closed ball in $\R^N$
centered at $x_0$ of radius $R$}\\
B_R & \hbox{the closed ball in $\R^N$
centered at $0$ of radius $R$}\\
\langle\,\cdot\;, \cdot\,\rangle & \hbox{the scalar product in
$\R^N$} \\
|\cdot| & \hbox{the Euclidean norm in $\R^N$}\\
\R_+& \hbox{the set of nonnegative real numbers}\\
\mathcal{B}(\R^N) & \hbox{the $\sigma$--algebra of Borel subsets
of $\R^N$}\\
\cchi_E & \hbox{the characteristic function of the set $E$}\\
\end{array}
\]

\vspace{1ex} Given a subset $U$ of $\R^N$, we denote by $\overline
U$ its closure. We furthermore say that $U$ is {\em compactly
contained} in a subset $V$ of $\R^N$ if $\overline U$ is compact
and contained in $V$. If $E$ is a Lebesgue measurable subset of
$\R^N$, we denote by $|E|$ its $N$--dimensional Lebesgue measure,
and qualify $E$ as {\em negligible} whenever $|E|=0$. We say that
a property holds {\em almost everywhere} ($a.e.$ for short) on
$\R^N$ if it holds up to a negligible set. We will write
$\varphi_n\ucv\varphi$ on $\R^N$ to mean that the sequence of
functions $(\varphi_n)_n$ uniformly converges to $\varphi$ on
compact subsets of $\R^N$.\par

With the term {\em curve}, without any further specification, we
refer to a Lipschitz--continuous function from some given interval
$[a,b]$ to $\R^N$. The space of all such curves is denoted by
$\D{Lip}([a,b],\R^N)$, while $\D{Lip}_{x,y}([a,b],\R^N)$ stands
for the family of curves $\gamma$  joining  $x$ to $y$, i.e. such
that $\gamma (a)=x$ and $\gamma (b)=y$, for any fixed $x$, $y$ in
$\R^N$. We denote by $W^{1,1}([a,b],\R^N)$ the space of absolutely
continuous curves defined in $[a,b]$. Given a curve $\gamma$
defined on some interval $[a,b]$, a curve $\gamma'$ defined on
$[a',b']$ will be called a {\em reparametrization of $\gamma$} if
there exists an order preserving Lipschitz--continuous map
$f:[a',b']\to[a,b]$ surjective and such that
$\gamma'=\gamma\,\comp\,f$. The Euclidean length of a curve
$\gamma$ is denoted by $\hh^1(\gamma)$.\par

For a measurable function $g:I\to\R^N$, $\|g\|_\infty$ stands for
$\sqrt{\sum_{i=0}^N \|g_i\|_{L^\infty(I)}^2}$, where $g_i$ and
$\|g_i\|_{L^\infty(I)}$ denotes the $i$--th component of $g$ and
the $\D L^\infty$--norm of $g_i$ respectively.
\par

%Unless otherwise specified, the term (sub, super) solution to some
%PDE equation is understood in the viscosity sense. Given a
%continuous function $g$ defined in $\R^k$ and $x_0 \in \R^k$, we
%denote by $D^+g(x_0)$ (resp. $D^-g(x_0)$) the superdifferential
%(resp. the subdifferential) of $g$ at $x_0$, i.e. the (possibly
%empty) set made up by the differentials of viscosity test function
%from above (resp. from below) of $g$ at $x_0$. Note that, in the
%case where $g$ is convex, $D^-g$ coincides with the usual
%subdifferential of convex analysis. When $g$ is defined on
%$\R^m\times\R^k$ and $(x_0,p_0)\in\R^m\times\R^k$, we will denote
%by $D^-_p g(x_0,p_0)$ the subdifferential of the function
%$g(x_0,\cdot)$ at $p_0$. For a function $g: \R^k \to (- \infty,
%+\infty]$, we denote by $\D{dom}(g)$ its effective domain, i.e.
%the subset of $\R^k$ where $g$ is finite valued.\par

\smallskip

Throughout the paper, $(\Omega,\F, \PP)$ will denote a {\em
probability space},  where $\PP$ is the probability measure  and
$\F$ the $\sigma$--algebra of $\PP$--measurable sets. A property
will be said to hold {\em almost surely} ($a.s.$ for short) in
$\omega$ if it holds up to a subset of probability 0. We will
indicate by $L^p(\Omega)$, $p\geq 1$, the usual Lebesgue space on
$\Omega$ with respect to $\PP$. If $f\in L^1(\Omega)$, we write
$\EE(f)$ for the mean of $f$ on $\Omega$, i.e. the quantity
$\int_\Omega f(\omega)\,\dd \PP(\omega)$.\par

We qualify  as {\em measurable} a map from $\Omega$ to itself,  or
 to a topological space $\M$ with Borel
$\sigma$--algebra $\mathcal{B}(\M)$,  if the inverse image of any
set in $\F$ or in $\mathcal{B}(\M)$  belongs to $\F$. This object
will be also called {\em random variable} with values in $\M$.
\par

We will be particulary interested in the case where the range of a
random variable is a {\em Polish}  space, namely  a complete and
separable metric space. By $\D C(\R^N)$ and
$\D{Lip}_\kappa(\R^n)$ we will denote the Polish
 space of continuous and
Lipschitz--continuous real functions, with Lipschitz constant less
than or equal to $\kappa>0$, defined in  $\R^N$, both endowed with
the metric of the uniform convergence on compact subsets
of $\R^N$. %, defined by
%\begin{equation}\label{def d}
%d(f,g):=\sum_{n=1}^\infty
%\frac{1}{2^n}\,\frac{\|f-g\|_{L^\infty(B_n)}}{\|f-g\|_{L^\infty(B_n)}+1}\qquad
%f,g\in\D C(\R^k).
%\end{equation}
We will use the expressions {\em continuous random function}, {\em
$\kappa$--Lipschitz random function}, respectively, for the
previously introduced random variables. Actually, we will usually  omit
$\kappa$ and  simply write {\em Lipschitz random function}. The
following characterization of random continuous functions holds,
see \cite{DavSic08}

\begin{prop}\label{prop misurability}
Let $\omega \mapsto v(\cdot,\omega)$ be a map from $\Omega$ to
$\D{C}(\R^N)$. The following are equivalent facts:
\begin{itemize}
    \item[\em (i)]
    $v$  is a random continuous function;\smallskip
    \item[\em (ii)] for every $x\in\R^N$, the map
    $\omega\mapsto v(x,\omega)$ is measurable in $\Omega$;\smallskip
    \item[\em (iii)] the map $(x,\omega)\mapsto v(x,\omega)$
    is jointly measurable
in $\R^N\times\Omega$, i.e. measurable  with respect to the
product $\sigma$--algebra $\mathcal B(\R^N)\otimes\F$.\smallskip
\end{itemize}
\end{prop}

Throughout the paper $(\tau_x)_{x\in\R^N}$ will denote a {\em
$N$--dimensional dynamical system}, defined as a family of
mappings $\tau_x:\Omega\to\Omega$ which satisfy the following
properties:
\begin{enumerate}
\item[{\em (1)}] the {\em group property:} $\tau_0=id$,\quad
$\tau_{x+y}=\tau_x\comp\tau_y$;

\item[{\em (2)}] the mappings $\tau_x:\Omega\to\Omega$ are
measurable and measure preserving, i.e. $\PP(\tau_x E)=\PP(E)$ for
every $E\in\F$;

\item[{\em (3)}] the map $(x,\omega)\mapsto \tau_x\omega$ from
$\R^N\times\Omega$ to $\Omega$ is jointly measurable, i.e.
measurable  with respect to the product $\sigma$--algebra
$\mathcal B (\R^N)\otimes\F$.
\end{enumerate}

We will moreover assume that $(\tau_x)_{x\in\R^N}$ is {\em
ergodic,} i.e. that one of the following equivalent conditions
hold:
\begin{itemize}
\item[{\em (i)}] every measurable function $f$ defined on $\Omega$
such that, for every $x\in\R^N$, $f(\tau_x\omega)=f(\omega)$ a.s.
in $\Omega$, is almost surely constant; \item[{\em (ii)}] every
set $A\in\F$ such that $\PP(\tau_x A\,\Delta\, A)=0$ for every
$x\in\R^N$ has probability either 0 or 1, where $\Delta$ stands
for the symmetric difference.
\end{itemize}

Notice that for  any vector subspace $V \subset \R^N$,
$(\tau_x)_{x \in V}$ is still a dynamical system  on $\Omega$, but
ergodicity does not hold in general.

Given   a random variable $f:\Omega\to\R$, for any fixed
$\omega\in\Omega$ the function $x\mapsto f(\tau_x\omega)$ is said
to be a {\em realization of $f$.} The following properties follow
from Fubini's Theorem, see \cite{JiKoOl}:  if $f\in L^p(\Omega)$,
then $\PP$--almost all its realizations belong to
$L^p_{loc}(\R^N)$; if $f_n \rightarrow f$ in $L^p(\Omega)$, then
$\PP$--almost all realizations of $f_n$ converge to the
corresponding realization of $f$ in $L^p_{loc}(\R^N)$. The
Lebesgue spaces on $\R^N$ are understood with respect to the
Lebesgue measure.\par

The next lemma guarantees that a modification of a random variable
on a set of zero probability does not affect its realizations on
sets of positive Lebesgue measure on $\R^N$, almost surely in
$\omega$. The proof is based on Fubini's Theorem again, see Lemma
7.1 in \cite{JiKoOl}.

\begin{lemma}\label{lemma utile}
Let $\widehat\Omega$ be a set of full measure in $\Omega$. Then
there exists a set of full measure
$\Omega'\subseteq\widehat\Omega$ such that for any
$\omega\in\Omega'$ we have $\tau_x\omega\in\widehat\Omega$ for
almost every $x\in\R^N$.
\end{lemma}

Next we state the Birkhoff Ergodic Theorem for ergodic
$N$--dimensional dynamical systems. It establishes a relation
between statistical and spatial means.

\begin{teorema}[\bf Birkhoff Ergodic Theorem] Let $(\Omega,\F,\PP)$ be a probability space and
$(\tau_x)_{x\in\R^N}$ a group of translations as above. Then, for
any $f\in L^1(\Omega)$, the limit
\[
f^*(\omega)= \lim_{t\to +\infty} \int_{t E}\!\taglio
\,f(\tau_x\omega)\,\dd x\qquad
\]
exists and is invariant with respect to $(\tau_x)_{x\in\R^N}$ a.s.
in $\omega$, where $E$ is any Borel subset of $\R^N$ with $|E|
>0$. Moreover $\EE(f^*)=\EE(f)$. If,
in addition, $(\tau_x)_{x\in\R^N}$ is ergodic, then
$f^*(\omega)=\EE(f)$ a.s. in $\omega$.
\end{teorema}

We will also need the following subadditive ergodic theorem.

\begin{teorema}[\bf Kingman's Subadditive Ergodic Theorem]
Let $\{f_{m,n}\,:\, 0\leq m\leq n\}$ be random variables which
satisfy the following properties:
\begin{itemize}
    \item[\em (a)] $f_{0,0}=0$ and $f_{m,n}\leq
f_{m,k}+f_{k,m}$ for every ${m\leq k\leq n;}$\medskip
    \item[\em (b)]  $\{f_{m,m+k}\,:\, m\geq 0,\,k\geq 0\}$
 have the same distribution law than $\{f_{m+1,m+k+1}\,:\, m\geq 0,\,k\geq
 0\}$, i.e. for every $0\leq m_1<\dots<m_n$,  $0\leq k_1<\dots<k_n$, $n\in\N$
 \[
\displaystyle{ \PP\left(\cap_{i=1}^{n}\,
f_{m_1,m_1+k_1}^{-1}(A_i)\right) = \PP\left(\cap_{i=1}^{n}
\,f_{m_1+1,m_1+k_1+1}^{-1}(A_i)\right)}
\]
for any open subset $A_i$ of $\R$;\medskip
    \item[\em (c)] $\int_\Omega \left(f_{0,1}(\omega)\right)^+\,\dd \PP(\omega)<+\infty.$\medskip
\end{itemize}
Then the following holds:
\begin{itemize}
    \item[\em (i)] $\displaystyle{\mu:=\lim_{n\to\infty}\frac{1}{n}\int_\Omega
f_{0,n}(\omega)\,\dd
\PP(\omega)=\inf_{n\in\N}\frac{1}{n}\int_\Omega
f_{0,n}(\omega)\,\dd \PP(\omega)}\in [-\infty,+\infty)$;\smallskip
    \item[\em (ii)]
    $\displaystyle{f_\infty(\omega):=\lim_{n\to\infty}\frac{f_{0,n}(\omega)}{n}}$ exists
    for $\PP$--almost every $\omega\in\Omega$;\smallskip
    \item[\em (iii)] $\displaystyle{\int_\Omega f_\infty(\omega)\,\dd
    \PP(\omega)=\mu}$ and, if $\mu>-\infty$, then
    \[
\frac{f_{0,n}}{n}\to f_\infty\quad\hbox{in $L^1(\Omega)$.}
    \]
\end{itemize}
\end{teorema}
\medskip
\end{section}

\begin{section}{Stationary Random Variables} \label{statvar}

In this section we recall the notion of stationarity for random
functions and random sets. These objects are of crucial relevance
for the extension of Lax--type formulae to the stationary ergodic
setting, see Propositions \ref{prop Lax in S} and \ref{prop Lax in
S0}. Then we will proceed to give the definition and to study the
properties of Lipschitz random functions with stationary
increments. Some of these results have been already proved in
\cite{DavSic08} for $N=1$. Their generalization to higher
dimensions is more subtle and requires additional tools, whose
presentation has been postponed to the Appendix, as well as those
proofs that on such tools are based, i.e. Theorems \ref{ascoli},
\ref{teorema mean 0} and \ref{teorema constant
mean}.\smallskip\par

%In order to keep this section as short as possible, we have
%postponed to the Appendix the presentation of this auxiliary
%material, as well as the proofs of those results that on such
%tools are based, specifically Theorems \ref{ascoli}, \ref{teorema
%mean 0} and \ref{teorema constant mean}.\smallskip\par

A jointly measurable function $v$ defined in $\R^N \times \Omega$
is said {\em stationary } if, for every $z \in \R^N$, there exists
a set $\Omega_z$ with probability $1$ such that for every
$\omega\in\Omega_z$
\[ v(\cdot + z, \omega)= v(\cdot, \tau_z \omega) \quad\text{on $\R^N$}
\]
It is clear that a real random
variable $\phi$  gives rise to a stationary function $v$
 by setting $v(x,\omega)=
\phi(\tau_x \omega)$. Conversely, according to Proposition 3.1 in
\cite{DavSic08}, a stationary function $v$ is, a.s. in $\omega$,
the realization of the measurable function $\omega \mapsto
v(0,\omega)$. More precisely,  there exists  a set $\Omega'$ of
probability $1$ such that for every $\omega\in\Omega'$
\begin{equation}\label{staz}
   v(x,\omega)=  v(0, \tau_x \omega) \quad\text{for a.e. $x\in\R^N$.}
\end{equation}

With the term {\em (graph--measurable ) random set} we indicate  a
set--valued function $X:\Omega\to\mathcal B(\R^N)$ with
\[
\Gamma(X):=\left\{(x,\omega)\in\R^N\times\Omega\,:\,x\in
X(\omega)\,\right\}
\]
jointly measurable in $\R^N\times\Omega$. A random set $X$ will be
qualified as {\em stationary} if for every   $z\in\R^N$, there
exists a set $\Omega_z$ of probability 1 such that
\begin{equation}\label{def stationary set}
X(\tau_z\omega)=X(\omega)-z\qquad\hbox{for every
$\omega\in\Omega_z$.}
\end{equation}

We use a stronger notion of measurability, which is usually  named
in the literature after Effros,  to define a {\em closed random
set}, say $X(\omega)$. Namely we require $X(\omega)$ to be a
closed subset of $\R^N$ for any $\omega$ and
\begin{equation*}\label{effros}
    \{\omega \, : \, X(\omega) \cap K \neq \emptyset \} \in \F
\end{equation*}
with $K$ varying  among the compact (equivalently, open) subsets
of $\R^N$. This condition can be analogously expressed by saying
that $X$ is measurable with respect to the Borel $\sigma$--algebra
related to the Fell topology on the family of closed subsets of
$\R^N$. This, in turn, coincides with the Effros
$\sigma$--algebra. If $X(\omega)$ is measurable in this sense then
it is also graph--measurable, see \cite{Molchanov} for more
details.

A closed random set $X$ is called stationary if it, in addition,
satisfies \eqref{def stationary set}. Note that in this event the
set $\{\omega\,:\,X(\omega)\not=\emptyset\,\}$, which is
measurable by the Effros measurability of $X$, is invariant with
respect to the group of translations $(\tau_x)_{x\in\R^N}$ by
stationarity, so it has probability either 0 or 1 by the
ergodicity assumption.

A convenient way to produce random closed (stationary) sets in
$\R^N$ is indicated by the next result, see \cite{DavSic08} for a
proof.

\begin{prop}\label{rando} Let $f$ be a continuous  random
function and $C$ a closed subset of $\R$. Then
\begin{equation*}\label{rando_1}
    X(\omega):= \{x \, :\, f(x,\omega)\in C\}
\end{equation*}
is a closed  random set in $\R^N$. If in addition $f$ is
stationary, then $X$ is stationary.
\end{prop}

%\begin{oss}\label{oss fat}
%We record for later use that, for a positive radius $r$, the set
%\[
%X_r(\omega):=X(\omega)+B_r,\quad\omega\in\Omega,
%\]
%is  closed random whenever $X$ is an almost surely nonempty closed
%random set, see Theorem III-9 in \cite{CaVa77} or Theorem 2.3 in
%\cite{Molchanov}.
%\end{oss}

For a random stationary set $X$ is immediate, by exploiting that
the maps $\{\tau_x\}_{x \in \R^N}$  are measure preserving, that
$\PP (X^{-1}(x))$ does not depend on $x$, where
\[X^{-1}(x)=\{\omega \, :\, x \in X(\omega)\}.\] Such quantity
 will be called {\em volume
fraction} of $X$ and denoted by $q_X$. Note that to any measurable
subset $\Omega'$ of $\Omega$ it can be associated a stationary set
$Y$  through the formula
\begin{equation*}%\label{stat1}
    Y(\omega):= \{x \,:\, \tau_x \omega \in \Omega'\}.
\end{equation*}
In this case $Y^{-1}(x)= \tau_{-x}\Omega'$, and so
$q_Y=\PP(\Omega')$.

The following classical result, which can be obtained as a direct
application of Fubini's theorem, will play a relevant role in what
follows, see \cite{Molchanov}.
\begin{teorema}[\bf Robbins' Theorem] Let $X$ be a  random set in $\R^N$.
If $\mu$ is a locally finite measure on Borel sets, then $\mu(X)$
is a random variable and
\begin{equation*}%\label{eq Robbins}
\int_\Omega \mu(X(\omega))\,\dd \PP=\int_{\R^N} \PP (X^{-1}(x))\,
\dd \mu,
\end{equation*}
in the sense that if one side is finite, then so is the other and
they are equal.
\end{teorema}

We next exploit the ergodicity assumption to get, through the
Birkhoff Ergodic Theorem,  an interesting  information on the
asymptotic structure of stationary sets, yielding in particular
that stationary sets are spread with some uniformity in the space.
We refer the reader to \cite{DavSic08} for the proof.

\begin{prop}\label{prop stationary random set}
Let $X$ be an almost surely nonempty closed stationary set in
$\R^N$. Then for every $\eps>0$ there exists $R_\eps>0$ such that
\[
\lim_{r\to +\infty} \frac{|\left(X(\omega)+B_R\right)\cap
B_r|}{|B_r|}\geq 1-\eps\qquad\hbox{a.s. in $\Omega$,}
\]
whenever $R\geq R_\eps$.
\end{prop}

\indent Given a Lipschitz random function $v$, we set
\[
\Delta_v(\omega):=\left\{\,x\in\R^N\,:\,v(\cdot,\omega)\ \hbox{is
differentiable at $x$}\, \right\}.
\]

\begin{definizione}\label{def v with st increments}
A random Lipschitz function $v$ is said to have
 {\em stationary increments} if, for every
    $z\in\R^N$, there exists a set $\Omega_z$ of probability 1
    such that
\begin{equation*}%\label{def stationary increments}
    v(x+z,\omega)-v(y+z,\omega)=v(x,\tau_z
    \omega)-v(y,\tau_z\omega)\quad\hbox{for all  $x,y\in\R^N$}
\end{equation*}
    for every $\omega\in\Omega_z$.
\end{definizione}

    The following holds:

\begin{prop}\label{prop Delta_v random}
Let $v$ be a Lipschitz random function, then $\Delta_v$ is a
random set. In addition, it is stationary with volume
 fraction $1$ whenever $v$ has stationary increments.
\end{prop}

\begin{dimo} The property of $\Delta_v$ of being a random set can
be proved via standard measure theoretic arguments, see for
instance Lemma 2.5 in \cite{FatFig} for a short proof. If $v$ has
stationary increments then, for any fixed $z\in\R^N$,
\[
v(\cdot+z,\omega)-v(\cdot,\tau_z\omega) \quad\hbox{is constant on
$\R^N$}
\]
whenever $\omega$ belongs to some set $\Omega_z$ with probability
$1$. This implies that $x+z$ is a differentiability point for
$v(\cdot,\omega)$ if and only if $x$ is a differentiability point
for $v(\cdot,\tau_z\omega)$, which, in turn,  means that
$\Delta_v$ is a stationary.  Since $\Delta_v(\omega)$ has full
measure in $\R^N$ for every $\omega$ by Rademacher's Theorem,
Robbin's Theorem with $\mu$ equal the Lebesgue measure restricted
to some ball of $\R^N$ readily implies that it volume fraction is
equal to 1.
\end{dimo}

\ \\
\indent Next, we state an important stability result for random
functions with stationary increments which are equiLipschitz, i.e.
that all take values in $\D{Lip}_\kappa(\R^N)$ for some fixed
$\kappa>0$. \medskip
\begin{teorema}\label{ascoli} Let $(v_n)_n$ be an
equiLipschitz   sequence of
 random functions with stationary increments.
Then there exist a random Lipschitz function $v$ with stationary
increments, a sequence $w_k= \sum_{n \geq n_k} \lambda_n^k v_n$ of
finite convex combinations of the $v_n$  and  a sequence $g_k$ of
real random variables  such that
\[ w_k(\cdot,\omega) + g_k(\omega) \ucv v(\cdot,\omega) \quad\text{a.s. in $\omega$.}\]
In addition the sequence of indices $(n_k)_k$ can be taken
diverging.\\
\end{teorema}

Let $v$ be a Lipschitz random function with stationary gradient.
For every fixed $x\in\R^N$, the random variable $Dv(x,\cdot)$ is
well defined on $\Delta_v^{-1}(x)$, which has probability 1 since
$\Delta_v$ is a stationary set with volume fraction 1.
Accordingly, we can define the mean $\EE(Dv(x,\cdot))$, which is
furthermore independent of $x$, see Proposition \ref{dif}--{\em
(i)}. In the sequel, we will be especially interested in the case
when this mean is zero.

\begin{definizione}\label{def ammissible}
A Lipschitz random function will be called {\em admissible} if it
has stationary increments and gradient with mean 0.
\end{definizione}

We state two characterizations of admissible random functions.

\begin{teorema}\label{teorema mean 0}
A Lipschitz random function $v$  with stationary increments has
gradient with vanishing mean if and only if it is almost surely
sublinear at infinity, namely
\begin{equation}\label{mean1}
    \lim_{|x|\to
    +\infty}\frac{v(x,\omega)}{|x|}=0\qquad\hbox{a.s. in $\omega$.}
\end{equation}
\end{teorema}
\smallskip
\begin{teorema}\label{teorema constant mean}
A Lipschitz random function  $v$ with stationary increments has
gradient  with vanishing mean if and only if
\begin{equation}\label{ciao}
   x\mapsto \EE(v(y,\cdot)- v(x,\cdot))=0 \quad\text{for any $x,y\in\R^N$.}\bigskip
\end{equation}
\end{teorema}

Finally, we show that any stationary Lipschitz random function is
admissible. Notice that in this case $\EE(v(x,\cdot))$ is
independent of $x$, so when such a quantity is finite this is just
a consequence of Theorem \ref{teorema constant mean}.

\begin{teorema}
Any stationary Lipschitz random function $v$ is admissible.
\end{teorema}

\begin{dimo}
Clearly, a stationary Lipschitz random function has stationary
increments. To prove the assertion, it is therefore enough to show
that $v(\cdot,\omega)$ is almost surely sublinear at infinity, in
view of Theorem \ref{teorema mean 0}. Let $\kappa$ be a Lipschitz
constant for $v(\cdot,\omega)$ for every $\omega$. The stationary
character of $v$ means, cf. \eqref{staz}, that for any $\omega$ in
a set  of probability 1
\[
v(x,\omega)=v(0,\tau_x\omega)\qquad\hbox{for a.e. $x\in\R^N$.}
\]
We claim that there exists a constant $M$ such that
$E:=\{\omega\in\Omega\,:\,|v(0,\omega)|\leq M\,\}$ has positive
probability. Indeed, if this were not the case, we would have that
$E_\infty:=\{\omega\,:\, |v(0,\omega)|=+\infty\,\}$ has
probability 1. An application of Robbins' Theorem with
$\mu=\leb^N\restr B_r$ for every fixed $r>0$ would imply that the
stationary random set
\[
F_\infty(\omega)=\{x\in\R^N\,:\,\tau_x\omega\in E_\infty,\,
v(x,\omega)=v(0,\tau_x\omega)\,\}
\]
is of full measure in $\R^N$ a.s. in $\omega$, yielding
$v(\cdot,\omega)\equiv +\infty$ a.s. in $\omega$, a contradiction.

Let us fix $M$ such that $E$ has positive probability. Then  the
closed stationary random set
\[
C(\omega)=\{x\in\R^N\,:\,|v(x,\omega)|\leq M\,\}
\]
is almost surely nonempty. Accordingly, by Proposition \ref{prop
stationary random set} there exists a set $\Omega_0$ probability 1
such that, for every $\omega\in\Omega_0$,
\begin{equation*}
\lim_{r\to +\infty} \frac{|\left(C(\omega)+B_n\right)\cap
B_r|}{|B_r|}> 1-\eps_n,
\end{equation*}
where $(\eps_n)_n$ is a sequence decreasing to 0. Fix
$\omega\in\Omega_0$. Then for every $x\in\R^N$ with $|x|$ large
enough, we have
\[
{|\left(C(\omega)+B_n\right)\cap B_{2|x|}|}>
(1-\eps_n)\,{|B_{2|x|}|}.
\]
For $n$ sufficiently large  $B_{2 |x|(\eps_n)^{1/N}}(x)\subseteq
B_{2|x|}$, and from the above inequality we infer
\[
B_{2 |x|(\eps_n)^{1/N}}(x)\cap(C(\omega)+B_n)\not=\emptyset,
\]
i.e. there exists $y=y(x,n)$ in $C(\omega)$ such that $|y-x|<2
|x|(\eps_n)^{1/N}+n$. Since $|v(y,\omega)|\leq M$, we get
\begin{eqnarray*}
|v(x,\omega)|\leq |v(x,\omega)-v(y,\omega)|+|v(y,\omega)|\leq
\kappa\left(2|x|(\eps_n)^{1/N}+n\right)+M.
\end{eqnarray*}
From this  we obtain
\[
\limsup_{|x|\to +\infty}\frac{|v(x,\omega)|}{|x|}\leq
2\kappa\,(\eps_n)^{1/N},
\]
and the claim follows letting $n\to +\infty$.
\end{dimo}
\ \\
\end{section}

\begin{section}{ Stochastic Hamilton--Jacobi equations}\label{sez HJ}

We consider an Hamiltonian
\[
H:\R^N\times\R^N\times\Omega\to\R
\]
satisfying the following conditions:

\begin{itemize}
    \item[(H1)] the map $\omega\mapsto H(\cdot,\cdot,\omega)$
    from $\Omega$ to the Polish space $C(\R^N\times\R^N)$ is
    measurable;\smallskip
    \item[(H2)] for every $(x,\omega)\in\R^N\times\Omega$,
    $\ H(x,\cdot,\omega)\ \hbox{is  convex on $\R^N$;}$ \smallskip
     \item[(H3)]  there
     exist two superlinear functions $\alpha,\beta:\R_+\to\R$ such
     that
     \[
     \alpha\left(|p|\right)\leq H(x,p,\omega)\leq \beta\left(|p|\right)\qquad\hbox{for all
     $(x,p,\omega)\in\R^N\times\R^N\times\Omega$;}
     \]
 %    \smallskip
      \item[(H4)] for every $(x,\omega)\in\R^N\times\Omega$,
      the set of minimizers of $H(x,\cdot,\omega)$
    has empty interior;\smallskip
    \item[(H5)] $H(\cdot+z,\cdot,\omega)=H(\cdot,\cdot,\tau_z\omega)$ for
    every $(z,\omega)\in\R^N\times\Omega$.
\end{itemize}

\begin{oss}
Condition (H3) is equivalent to saying that $H$ is superlinear and
locally bounded in $p$, uniformly with respect to $(x,\omega)$. We
deduce from (H2)
\begin{equation}\label{lippo}
  | H(x,p,\omega)-H(x,q,\omega)| \leq L_R |p-q| \quad\text{for all
  $x$, $\omega$, and $p$, $q$ in $B_R$},
\end{equation}
where $L_R=\sup\{\, |H(x,p,\omega)|\,:\,(x.\omega) \in \R^N \times
\Omega, \,|p|\leq R+2\,\},$ which is finite thanks to (H3). For a
comment on hypothesis (H4), see Remark \ref{oss h4}.
\end{oss}

\begin{oss}\label{oss almost-periodic}
Any given periodic, quasi--periodic or almost--periodic
Hamiltonian $H_0:\R^N\times\R^N\to\R$ can be seen as a specific
realization of a suitably defined stationary ergodic Hamiltonian,
cf. Remark 4.2 in \cite{DavSic08}.
\end{oss}

For every $a\in\R$, we are interested  in the stochastic
Hamilton--Jacobi equation
\begin{equation}\label{eq HJa}
H(x,Dv(x,\omega),\omega)=a\qquad\hbox{in $\R^N$.}
\end{equation}
The analysis performed on it in \cite{DavSic08} stays valid in the
present multidimensional setting, with minor adjustments. We
basically refer to it,  just recalling the main  items and
pointing out the main differences. \par

We say that a Lipschitz random function is a {\em solution} (resp.
{\em subsolution}) of \abra{eq HJa} if it is a viscosity solution
(resp.  a.e. subsolution) a.s. in $\omega$ (see \cite{BCD97, Ba94}
for the definition of viscosity (sub)solution in the deterministic
case). Notice that any such subsolution takes value in
$\D{Lip}_{\kappa_a}(\R^n)$, where
\begin{equation*}%\label{def kappa_a}
\kappa_a:=\sup\{\,|p|\,:\,H(x,p,\omega)\leq a\ \hbox{for some
$(x,\omega)\in\R^N\times\Omega$}\,\},
\end{equation*}
and this quantity  is finite thanks to (H3). We are interested in
the class of {\em admissible subsolutions}, hereafter denoted by
$\S_a$, i.e. random functions  with stationary increments and zero
mean gradient that are subsolutions of \eqref{eq HJa}. An
admissible solution will be also named {\em exact corrector},
remembering its role in homogenization. Further, for any $\delta
>0$, a random function $v_\delta$ will be called a {\em
$\delta$--approximate corrector} for the equation \eqref{eq HJa}
if it belongs to $\S_{a+\delta}$ and satisfies the inequalities
\begin{equation}\label{approxima}
    a - \delta \leq H(x, Dv_\delta(x,\omega),\omega) \leq a +\delta
\end{equation}
in the viscosity sense a.s. in $\omega$. We say that \abra{eq HJa}
has {\em approximate correctors} if it admits
$\delta$--approximate correctors for any $\delta >0$. \par

The following stability property of  admissible subsolutions is a
consequence of Theorem \ref{ascoli}  along with the remark that,
if in the convergence established there the approximating random
functions are admissible, the limit too keeps this property. The
proof is analogous to that of Theorem 4.3 in \cite{DavSic08}.

\begin{teorema}\label{stabile} Let $(a_n)_{n}$ be a sequence of real numbers and
$v_n$ a random function in $\S_{a_n}$ for each $n$. If $a_n$
converges to some $a$, there exist $v \in\S_a$ and a sequence
$(w_k)_k$ made up by finite convex combinations of the $v_n$, up
to an additive real random variable, such that
\begin{equation*}\label{stabile1}
    w_k(\cdot,\omega) \ucv v(\cdot,\omega) \quad\text{a.s. in $\omega$.}
\end{equation*}
\end{teorema}

\medskip

We proceed by defining the {\em free} and the {\em stationary
critical value}, denoted by $c_f(\omega)$ and $c$ respectively, as
follows:
\begin{eqnarray}
 c &=& \inf\{a\in\R\,:\,\S_a\not =\emptyset\,\},  \label{def c}\\
 c_{f}(\omega)&=& \inf\left\{ a\in\R\,:\, \text{\abra{eq HJa} has a subsolution
$v\in\D{Lip}(\R^N)$} \right\}. \label{def cf}
\end{eqnarray}
We emphasize that in definition \abra{def cf} we are considering
{\em deterministic} a.e. subsolutions $v$ of the equation \abra{eq
HJa}, where $\omega$ is treated as a fixed parameter. The set
appearing at the right--hand side of \eqref{def c} is non void,
since it contains the value $\sup_{(x,\omega)} H(x,0,\omega)$,
which is finite thanks to (H3). Furthermore, notice that
$c_f(\tau_z\omega)=c_f (\omega)$ for every
$(z,\omega)\in\R^N\times\Omega$, so that, by ergodicity, the
random variable $c_f(\omega)$ is almost surely equal to a
constant, still denoted by  $c_f$. Hereafter we will write
$\Omega_f$ for the set of probability $1$ where $c_f(\omega)$ is
equal to  $c_f$. It is apparent that $c \geq c_f$.
%We will show in
%the sequel that the critical value can be equivalently defined
%through
%\begin{equation}\label{def cbis}
%   c= \inf\{ a \in \R \, :  \,\text{ \abra{eq HJa}  admits stationary
%   subsolutions}\}.
%\end{equation} i

In what follows, we mostly focus our attention on the {\em
critical equation}
\begin{equation}\label{eq critica}
H(x,Dv(x,\omega),\omega)=c\qquad\hbox{in $\R^N$.}
\end{equation}
It follows from Theorem \ref{stabile} that it admits admissible
subsolutions, i.e. $\S_c \neq \emptyset$. The relevance of the
critical value $c$ is given by the following result, see Theorem
4.5 in \cite{DavSic08} for the proof.

\begin{teorema}\label{correctors} The critical equation \eqref{eq
critica} is the unique among equations \eqref{eq HJa} either an
exact corrector or approximate correctors may exist. \end{teorema}

\smallskip

We introduce an intrinsic
 path distance, assuming   in next formulae
 $a \geq c_f$ and  $\omega \in \Omega_f$. We  start by defining
the sublevels
\[
Z_a(x,\omega):=\{p\,:\,H(x,p,\omega)\leq a\,\},
\]
and the related support functions  $\sigma_a(x,q,\omega)$ by
    \begin{equation*}\label{sigma}
    \sigma_a(x,q,\omega):=\sup\left\{\langle q,p\rangle\,:\,p \in Z_a(x,\omega)\,\right\}.
    \end{equation*}
It comes from \abra{lippo} (cf. Lemma 4.6 in \cite{DavSic08})
that, given $b>a$, we can find $\delta=\delta(b,a)>0$ with
\begin{equation}\label{spessore}
Z_a(x,\omega)+B_\delta\subseteq Z_b(x,\omega) \qquad\hbox{for
every $(x,\omega)\in\R^N\times\Omega_f$.}
\end{equation}
This property is used in the proof of Theorems \ref{teo metric
homogenization} and \ref{teo relation effective H}. It is also
needed in the proof of Theorem \ref{correctors}.  It is
straightforward to check that $\sigma_a$ is convex in $q$, upper
semicontinuous in $x$ and, in addition, continuous whenever
$Z_a(x,\omega)$ has nonempty interior or reduces to a point. We
extend the definition of $\sigma_a$ to
$\R^N\times\R^N\times\Omega$ by setting
$\sigma_a(\cdot,\cdot,\omega)\equiv 0$ for every
$\omega\in\Omega\setminus\Omega_f$. With this choice, the function
$\sigma_a$ is jointly measurable in $\R^N \times\R^N\times \Omega$
and enjoys the stationarity property
\[
\sigma_a(\cdot+z,\cdot,\omega)=\sigma_a(\cdot,\cdot,\tau_z\omega)\quad\hbox{for
every $z\in\R^N$ and $\omega\in\Omega$.}
\]
 Next, for every $a\geq c_f$,
we define the semidistance $S_a$ as
    \begin{equation}\label{eq S}
    S_a(x,y,\omega)=\inf\left\{\int_0^1 \sigma_a(\gamma(s),\dot\gamma(s),\omega)\,\dd
    s\,:\, \gamma\in\D{Lip}_{x,y}([0,1],\R^N)\, \right\}.
    \end{equation}
The function $S_a$ is measurable on $\R^N\times\R^N\times\Omega$
with respect to the product $\sigma$--algebra $\mathcal
B(\R^N)\otimes\mathcal B(\R^N)\otimes\F$, and satisfies the
following properties:
\begin{eqnarray*}
S_a(x,y,\tau_z\omega)&=&S(x+z,y+z,\omega)\\
S_a(x,y,\omega)&\leq& S_a(x,z,\omega)+S_a(z,y,\omega)\\
S_a(x,y,\omega)&\leq& \kappa_a|x-y|
\end{eqnarray*}
for all $x,y,z\in\R^N$ and $\omega\in\Omega$. According to
Proposition \ref{rando}, $S_a$ is a random (semi)distance, i.e. a
random variable taking values in the space of semidistances
endowed with the local uniform convergence in $\R^N \times \R^N$.
 We have (see \cite{FaSic03}):

\begin{prop}\label{prop S} Let $ a \geq c_f$ and $\omega\in\Omega_f$.
\begin{itemize}
    \item[\em(i)] For any $y\in\R^N$, the functions $S_a(y,\cdot,\omega)$ and
    $-S_a(\cdot,y,\omega)$ are both subsolutions of \abra{eq HJa}.
    \item[\em(ii)] A continuous function $\phi$ is a subsolution of \abra{eq HJa} if and
    only if
\[
\phi(x)-\phi(y)\leq S_a(y,x,\omega)\qquad\hbox{for all
$x,y\in\R^N$.}
\]

\end{itemize}
\end{prop}

\medskip
An immediate consequence of  the previous item (ii)  is that for
any cycle $\gamma$ defined in $[0,1]$ $\int_0^1
\sigma_a(\gamma(s),\dot\gamma(s),\omega)\,\dd
    s \geq 0$, whenever $a \geq c_f$.  We define  for every $\omega\in\Omega$
    the {\em classical (projected) Aubry set}  (cf.
\cite{FaSic03}), which  plays a  special role in the study of
equation \eqref{eq HJa} with $a=c_f$, as the collection of points
$y\in\R^N$ such that there is a sequence of cycles $\gamma_n$,
defined in $[0,1]$ and  based on $y$, with
\begin{equation*}\label{eq A}
\inf_n \int_0^1\sigma_{c_f}(\gamma_n,\dot\gamma_n,\omega)\,\dd s
=0 \quad\text{and} \quad \inf_n\,\hh^1(\gamma_n) >0
\end{equation*}
or, equivalently (cf. \cite[Lemma 5.1]{FaSic03}),
\begin{equation*}\label{eq A2}
\inf\left\{\int_0^1\sigma_{c_f}(\gamma,\dot\gamma,\omega)\,\dd
s\,:\,\gamma\in\D{Lip}_{y,y}([0,1],\R^N),\,\hh^1(\gamma)\geq\delta
\right\}=0\quad\hbox{for any $\delta>0$.}
\end{equation*}
 Hereafter we
will denote by $\A_{f}(\omega)$ the collection of points $y$ of
$\R^N$ enjoying one of the two equivalent conditions above. It is
 closed for every $\omega\in\Omega$. Given $\omega \in \Omega_f$, $ a \geq c_f$, and
 a closed subset $C$ of
 $\R^N$, a standard way for producing a subsolution of (\ref{eq HJa}) is by means of the following
 Lax formula
 \begin{equation}\label{lax}
    \inf\{w_0(y)+ S_a(y,x,\omega)\,:\,y \in C\},
\end{equation}
where $w_0$ is a function defined on $C$ which is 1
Lipschitz--continuous with respect to $S_a(\cdot,\cdot,\omega)$,
i.e.
\[
w_0(x)-w_0(y)\leq S_{a}(y,x,\omega)\qquad\hbox{for every $x,y\in
C$.}
\]
We recall that the function given above is also the maximal
subsolution taking the value $w_0$ on $C$ and hence a solution in
$\R^N\setminus C$. Furthermore we have (see \cite{FaSic03}):

\begin{teorema}\label{teo 6.7}
Let $\omega\in\Omega_f$. Then
\begin{itemize}
\item[\em (i)]  If
 $C \subset \A_{f}(\omega)$ then \eqref{lax}, with $a=c_f$ and $w_0$  $1$--Lipschitz continuous
   with respect to $S_{c_f}$,
   yields a solution on the whole $\R^N$.\medskip
\item[\em (ii)] If $U$ is a bounded open subset of $\R^N$, $a >
c_f$ and $w_0$ a function defined on $\partial U$ which is
$1$--Lipschitz continuous with respect to $S_{c_f}$, then
\eqref{lax} with $C$ replaced by $\partial U$ is the unique
viscosity solution of the Dirichlet Problem:
\begin{eqnarray*}
\begin{cases}
H(x,D \phi(x),\omega)=a&\qquad\hbox{in $U$}\\
\phi(x)=w_0(x)&\qquad\hbox{on $\partial U$.}
\end{cases}
\end{eqnarray*}

\item[\em (iii)] If $U$ is as above, $a = c_f$ and $w_0$ a
function defined on $\partial U \cup (U \cap \A_f)$ which is
$1$--Lipschitz continuous with respect to $S_{c_f}$, then
\eqref{lax} with $C$ replaced
   by $\partial U \cup (U \cap \A_f)$
is the unique  viscosity solution of the Dirichlet Problem:
\begin{eqnarray*}
\begin{cases}
H(x,D \phi(x),\omega)=c_f&\qquad\hbox{in $U \setminus \A_f $}\\
\phi(x)=w_0(x)&\qquad\hbox{on $\partial U \cup (U \cap \A_f)$.}
\end{cases}
\end{eqnarray*}
\end{itemize}
\end{teorema}

\bigskip
\indent We define for every $\omega\in\Omega$ the set of {\em
equilibria} as follows:
\begin{equation*}%\label{def equilibria}
\Eq(\omega):=\{y\in\R\,:\,\min_p H(y,p,\omega)=c_f\,\}.
\end{equation*}
The set $\Eq(\omega)$ is a (possibly empty)  closed subset of
$\A_{f}(\omega)$ (cf. \cite[Lemma 5.2]{FaSic03}).  It is apparent
that  ${c_f} \geq \sup_{x \in \R^N} \min_{p \in \R^N}
H(x,p,\omega)$ a.s. in $\omega$;  we point out that $\Eq(\omega)$
is nonempty if and only if  the previous formula holds with an
equality and the sup is a maximum. In this case, $\Eq(\omega)$ is
made up by the points  where such a maximum is attained. Note that
$\omega \mapsto  \sup_{x \in \R^N} \min_{p \in \R^N}
H(x,p,\omega)$ is a random variable and consequently, by
ergodicity,  almost surely constant.

\begin{oss}\label{oss h4}
The inclusion $\Eq(\omega)\subseteq \A_f(\omega)$ depends on the
fact that the $c_f$--sublevel $\{p\,:\,H(y,p,\omega)\leq {c_f}\}$
is non--void and has empty interior when $y \in \Eq(\omega)$. The
latter is a consequence of (H4), and this is actually the unique
point where such condition is used.
\end{oss}

We proceed giving a stochastic version of Lax formula in order to
recover the previous properties in our setting. Let $C(\omega)$ be
an almost surely nonempty stationary closed random set in $\R^N$.
Take  a Lipschitz random function $g$ and set, for $a\geq c_f$,
\begin{equation}\label{def u}
u(x,\omega):= \inf\{g(y,\omega)+S_a(y,x,\omega)\,:\,y\in
C(\omega)\,\}\quad \hbox{$x\in\R^N$,}
\end{equation}
where we agree that $u(\cdot,\omega)\equiv 0$ when either
$C(\omega)=\emptyset$ or the infimum above is $-\infty$. The
following holds:

\begin{prop}\label{prop Lax in S}
Let $g$ be a stationary Lipschitz random function and $C(\omega)$,
$u$ as  above. Let us assume that, for some $a\geq c_f$, the
infimum in \eqref{def u} is finite a.s. in $\omega$. Then $u$ is a
stationary random variable belonging to $\S_a$ and satisfies
$u(\cdot,\omega)\leq g(\cdot,\omega)$ on $C(\omega)$ a.s. in
$\omega$. Moreover, $u$ is a viscosity solution of \eqref{eq HJa}
in $\R^N\setminus {C(\omega)}$ a.s. in $\omega$.
\end{prop}

\smallskip When $g$ is itself an admissible subsolution of
\eqref{eq HJa}, we can state a stronger version  of the previous
result.

\begin{prop}\label{prop Lax in S0}
Let $g$ be a random function belonging to $\S_a$ and $C(\omega)$,
$u$ as above. Then $u$ belongs to $\S_a$. In addition, it is a
viscosity solution of \eqref{eq HJa} in $\R^N\setminus
{C(\omega)}$, and takes the value $g(\cdot,\omega)$ on $C(\omega)$
a.s. in $\omega$.
\end{prop}

As already pointed in the Introduction, the property of being $C$
a closed stationary set is of crucial importance to show that
formula \eqref{def u} defines an admissible Lipschitz random
function. The proofs of the above results are analogous to those
of Proposition 5.2. and 5.3 in \cite{DavSic08}, respectively,
where the case $N=1$ is considered.\smallskip

Later on, we will make use of the Lax--type formula \eqref{def u}
when the random source set is either  $\A_f$ or  $\Eq$. In order
to exploit the previous results, we will need the following

\begin{prop} \label{randomclosed}
$\Eq$ and $\A_f$ are closed random stationary sets.
\end{prop}

\begin{dimo}
For every $(x,\omega)\in\R^N\times\Omega$, let us set
$h(x,\omega)=\min_p H(x,p,\omega)$ and
\[
f(x,\omega)=\inf
    \left\{
    \int_0^1\sigma_{c_f}(\gamma,\dot\gamma,\omega)\,\dd
    s\,:\,\gamma\in\D{Lip}_{x,x}([0,1],\R^N),\,\hh^1(\gamma)\geq 1
    \right\}.
\]
Then $\Eq(\omega)=\{\,y\in\R^N\,:\,h(y,\omega)=c_f\,\}$ and
$\A_f(\omega)=\{\,y\in\R^N\,:\,f(y,\omega)=0\,\}$. In view of
Proposition \ref{rando}, it suffices to show that the stationary
functions $h$ and $f$ are jointly measurable in $(x,\omega)$ and
continuous in $x$ for every fixed $\omega$.

The continuity of $h$ in $x$ can be directly deduced from its very
definition by making use of assumptions (H1) and (H3); for the
measurability issue, simply notice that
$h(x,\omega)=\inf_{p_k\in\Q} H(x,p_k,\omega)$.

Let us consider the random variable $f$. It is easy to derive from
its definition that $f(x,\omega)-f(y,\omega)\leq
2\,\kappa_{c_f}|x-y|$ for every $x,\,y\in\R^N$, thus proving the
continuous character of $f(\cdot,\omega)$ for every fixed
$\omega\in\Omega$. For the measurability issue, it suffices to
show that the map $\omega\mapsto f(x,\omega)$ is measurable in
$\Omega$ for every fixed $x\in\R^N$ by Proposition \ref{prop
misurability}. To this purpose, let us fix $x\in\R^N$, and
consider the countable family $(\gamma_n)_n$ of polygonal loops
with vertexes in $x+\Q^N$, having $x$ as base point, and of
Euclidean length greater than 1. We claim that
\[
f(x,\omega)=\inf_{n}\int_0^1\sigma_{c_f}(\gamma_n,\dot\gamma_n,\omega)\,\dd
s,
\]
which clearly implies the asserted measurability of $f$. Indeed,
for any loop $\gamma$ having $x$ as base point and Euclidean
length greater than or equal to 1, an Euler--type approximation provides  a subsequence
$(\gamma_{n_k})_k$ of the family uniformly converging to $\gamma$
and such that
$
\sup_k\|\dot\gamma_{n_k}\|_\infty<+\infty.
$
In particular, the curves $\gamma_{n_k}$ weakly converge to
$\gamma$ in $W^{1,1}([0,1],\R^N)$. Being
$\sigma_{c_f}(\cdot,\cdot,\omega)$ upper semicontinuous for every
fixed $\omega$, by classical results of the Calculus of Variations
\cite{BuGiHi98} we derive
\[
\limsup_{k\to
+\infty}\int_0^1\sigma_{c_f}(\gamma_{n_k},\,\dot\gamma_{n_k},\omega)\,\dd
s
    \leq
\int_0^1\sigma_{c_f}(\gamma,\,\dot\gamma,\omega)\,\dd s.
\]
That yields
\[
f(x,\omega)
    \leq
\inf_{n}\int_0^1\sigma_{c_f}(\gamma_n,\dot\gamma_n,\omega)\,\dd s
    \leq
\int_0^1\sigma_{c_f}(\gamma,\,\dot\gamma,\omega)\,\dd
s\qquad\hbox{for every $\omega$},
\]
and the claim follows by taking the infimum of the right--hand
side term of the above inequality when $\gamma$ varies in the
family of loops with base point $x$ and of Euclidean length
greater than or equal to 1.
\end{dimo}
\ \\

\end{section}

\begin{section}{Stable norms} \label{stable}
In this section, we show the existence of  asymptotic norm--type
functions associated with $S_a$, whenever $a \geq c_f$, and
explore their link with the effective Hamiltonian $\overline
H$. Given $\eps>0$, we define
\begin{eqnarray*}
  S_a^\eps(x,y,\omega)=\inf\left\{\int_0^1\sigma_a(\gamma(t)/\eps,\dot\gamma(t),\omega)\,\dd
t\,:\,\gamma\in\D{Lip}_{x,y}([0,1],\R^N)\,\right\}
\end{eqnarray*}
for every $x,y\in\R^N$ and $\omega\in\Omega_f$, where we agree
that $S_a^\eps(\cdot,\cdot,\omega)\equiv 0$ when
$\omega\in\Omega\setminus\Omega_f$. Note that
$S_a^\eps(x,y,\omega)=\eps\,S_a(x/\eps,y/\eps,\omega)$.

\begin{teorema}\label{teo metric homogenization}
Let $a\geq c_f$. There exists a convex and positively
1--homogeneous function $\phi_a:\R^N\to\R$  such that
\begin{equation}\label{eq cv metrics}
S_a^\eps(x,y,\omega)\quad\underset{\eps\to
0}{\ucv}\quad\phi_a(y-x),\qquad\hbox{$x,y\in\R^N$}.
\end{equation}
for any $\omega$ in a set $\Omega_a$ of probability 1. In
addition,  $\phi_a$  is nonnegative for $a= c$, and nondegenerate,
i.e. satisfying  $\phi_a(\cdot)\geq \delta_a |\cdot|$ for some
$\delta_a>0$, when $a>c$.
\end{teorema}

With some abuse of terminology, we will refer to the function
$\phi_a$ appearing in the statement above as the {\em stable norm}
associated with $S_a$, in analogy with the case of periodic
Riemannian metrics. The above theorem states that $\phi_a$ is a
{\em Minkowski norm} (i.e. a norm which fails to be symmetric)
when $a>c$; it can possibly degenerate when $a=c$, in the sense
that $\phi_c$ may be identically 0 along some directions.\\

\begin{dimo}
The proof is basically divided in two parts. In the first half, we
essentially follow the arguments of \cite{ReTa00}, to which we
refer for the details (cf. also \cite{Souga99}). The second half,
based on a combined use of Egoroff's and Birkhoff Ergodic
Theorems, follows an argument provided in \cite{Souga-personal},
which is also needed in \cite{Souga99} to complete the proof of
Theorem 1.\smallskip\\
\indent Since for every $\omega$ the functions
$\left\{S^\eps_a(\cdot,\cdot,\omega)\right\}_{\eps>0}$ are
equiLipschitz--continuous, the local uniformity of the asserted
convergence  is a consequence of Ascoli--Arzel\'a Theorem, once we
show that there is pointwise convergence. \par

The first step  is to consider  the sequence of random variables $
S_a(0,n\,q,\omega)$, where $q$ is any vector of $\R^N$. The
subadditive decomposition through the double indexed random
variables $ S_a(m\,q,n\,q,\omega)$, $0\leq m\leq n$, allows to
apply the Subadditive Ergodic Theorem and  deduce the existence of
\[ \lim_{n\to\infty} \frac{1}{n}\,S_a(0,n\,q,\omega)=
\lim_{n\to\infty} \,S^{1/n}_a(0,q,\omega)\] for $\omega$ belonging
to some set $\Omega_q$ of probability $1$. The estimate $ |
S_a(0,nq,\omega)|\leq \kappa_a\,n|q|$, which holds for every
$\omega$, implies that such limit  is  almost surely finite. Since
for every fixed $\omega$ the functions
$\left\{S^\eps_a(\cdot,\cdot,\omega)\right\}_{\eps>0}$ are
equiLipschitz--continuous,   we  derive that the same limit is
attained by $S_a^{\eps}(0,q,\omega)$, as $\eps$ goes to $0$, for
$\omega \in \Omega_q$, and  stays unaffected passing from $\omega$
to $\tau_z\omega$ for all $z \in \R^N$, which in turn implies that
it is almost surely constant by ergodicity. By possibly redefining
$\Omega_q$ if necessary, we set
\begin{equation}\label{eq prima cv}
\phi_a(q)= \lim_{\eps\to 0} S_a^{\eps}(0,q,\omega),
\qquad\hbox{$\omega\in\Omega_q$.}
\end{equation}
 Then $\phi_a$ is Lipschitz--continuous with Lipschitz constant $\kappa_a$.
 By taking  a sequence $(q_n)_n$   dense in
$\R^N$ and exploiting  the equiLipschitz--continuity of
$\{S^\eps_a(0,\cdot,\omega)\}_{\eps>0}$ and $\phi_a(\cdot)$, we
see that the convergence in \abra{eq prima cv} takes place for any
$q \in \R^N$ whenever $\omega \in\widehat\Omega_a:=\cap_n
\Omega_{q_n}$. In addition, by the Subadditive ergodic Theorem
\[ S^\eps_a(0,q,\omega) \rightarrow \phi_a(q) \quad\hbox{in
$L^1(\Omega)$\qquad for any $q\in\R^N$.}\] Let us now  fix $x$,
$y$ in $\R^N$. Since
$S^\eps_a(x,y,\omega)=S^\eps_a(0,y-x,\tau_{x/\eps}\omega)$ a.s. in
$\omega$, and $\tau_{x/\eps}$ is measure preserving, we deduce
\[
\EE  (|S^\eps_a(0,y-x,\cdot)-\phi_a(y-x)|) = \EE
(|S^\eps_a(x,y,\cdot)-\phi_a(y-x)|),
\]
and so
\[
S^\eps_a(x,y ,\omega) \rightarrow \phi_a(y-x) \qquad\hbox{in
$L^1(\Omega)$.}
\]
%This shows that if we apply the Subadditive
%Ergodic Theorem to the family of random variables
%$\{S^\eps_a(x,y,\omega)\}_{\eps
%>0}$, we  get
%\[ \lim_{\eps \rightarrow 0} S^\eps_a(x,y,\omega)= \phi_a(y-x) \quad\text{a.s. in $\omega$.}\]

We now proceed to show that there exists a fixed set of
probability 1 on which this convergence also holds pointwise, for
any pair $x$, $y$ in $\R^N$. For this, we make a combined use of
Egoroff's and Birkhoff Ergodic Theorem.\smallskip

 Since the functions $q \mapsto S^\eps_a(0,q,\omega)$ are
equiLipschitz--continuous and locally equibounded for every
$\omega$  we deduce that, for every $r>0$,
\[
\sup_{|q|\leq 2r}
\left|S^\eps_a(0,q,\omega)-\phi_a(q)\right|\underset{\eps\to
0}{\tto} 0\qquad \omega\in\widehat\Omega_a.
\]
We use Egoroff's Theorem to make this convergence uniform in
$\omega$ on large sets, when $r \in \Q^+$ (the set of positive
rational numbers): for every $\delta>0$, we find a set $A_\delta$
with $\PP(\Omega\setminus A_\delta)\leq \delta$ such that
\begin{equation*}
\sup_{\omega\in A_\delta}\left(\sup_{|q|\leq 2r}
\left|S^\eps_a(0,q,\omega)-\phi_a(q)\right|\right)\underset{\eps\to
0}{\tto} 0,
\end{equation*}
for every $r \in \Q^+$. The Birkhoff Ergodic Theorem applied to
the function $\cchi_{A_\delta}$ yields the existence of a set
$\Omega^\delta$ of probability 1 with
\begin{equation*}%\label{eq cchi A_delta}
\lim_{R\to +\infty} \int_{B_R}\tagli
\cchi_{A_\delta}(\tau_x\omega)\,\dd x=\PP(A_\delta)\qquad\hbox{for
every $\omega\in\Omega^\delta$},
\end{equation*}
in other terms  for every $\omega \in \Omega^\delta$
\begin{equation}\label{eq cchi A_delta bis}
    \frac{\left|\{x\in\R^N\,:\,\tau_x\omega\in A_\delta\,\}\cap
B_{R}\right|}{|B_{R}|}
    \geq
    \PP(A_\delta)-\delta\geq 1-2\delta,
\end{equation}
whenever $R$ is large enough.  We set
$\displaystyle{\Omega_a:=\cap_{\delta\in\Q^+}\Omega^\delta}$.
Given
 $\omega_0 \in \Omega_a$,  for any $\alpha >0$ we can therefore find, according to
\abra{eq cchi A_delta bis}, a pair of positive numbers
$\delta(\alpha)$ and $R(\alpha)$ such that,  if $|z_0| \geq
R(\alpha)$, any ball centered at $z_0$ with radius exceeding
$\alpha |z_0|$ must intersect $\{x\,:\,\tau_x\omega_0\in
A_{\delta(\alpha)}\,\}$, or equivalently
\begin{equation}\label{tris}
    |z_0- z|\leq \alpha |z_0| \quad\text{for some $z$ with
$\tau_z \omega_0 \in A_{\delta(\alpha)}$.}
\end{equation}
Now fix $\alpha>0$, and pick a pair of points $x$,\,$y$ in $\R^N$.
We assume that they both belong to $B_r$ for some  $r \in \Q^+$.
Let $\eps_0$ be such that
\[
\sup_{\omega\in A_{\delta(\alpha)}}\left(\sup_{|q|\leq 2r}
\left|S^\eps_a(0,q,\omega)-\phi_a(q)\right|\right)\leq \alpha
\quad\text{for $\eps \leq \eps_0$}
\]
and \[ \frac{|x|}{\eps_0} > R(\alpha).\] We denote,  for $\eps
\leq \eps_0$, by $z_\eps$ a point such that \abra{tris} holds true
with $z_\eps$, $\frac{x}{\eps}$ in place of $z$, $z_0$,
respectively. Accordingly $|x - \eps z_\eps| \leq \alpha \, r$,
and for $\eps \leq \eps_0$ we have
\begin{eqnarray*}
|S^\eps_a(x,y,\omega_0)\!\!\!&-&\!\!\!\phi_a(y-x)|\leq
    |S^\eps_a(x,y,\omega_0)-S^\eps_a(\eps z_\eps,y,\omega_0)|\qquad\qquad\qquad\qquad\qquad\qquad\qquad\\
&+&|S^\eps_a(\eps z_\eps,y,\omega_0)-\phi_a(y-\eps z_\eps)|+
    |\phi_a(y-\eps z_\eps)-\phi_a(y-x)|  \\
    &\leq&
    2\,\kappa_a \, \alpha\,r
    +
    |S^\eps_a(0,y-\eps z_\eps,\tau_{z_\eps}\omega_0)-\phi_a(y-\eps z_\eps)|
    \leq
    \alpha\, (2 \kappa_a \,r+1).
\end{eqnarray*}
As $\alpha$ was arbitrarily chosen, we conclude that
\[
\lim_{\eps\to 0} S^\eps_a(x,y,\omega_0)=\phi_a(y-x),
\]
as desired.

It comes from its very definition that $\phi_a$ is positively
homogeneous. To prove that it is convex, we pick $\omega \in
\Omega_a$, $\lambda\in (0,1)$, $x$,\,$y$ in $\R^N$, and we pass to
the limit, as $\eps$ goes to $0$, in the inequality
\[
    S^\eps_a(0,\lambda\,x+(1-\lambda)\,y,\omega)
    \leq
    S^\eps_a(0,\lambda\,x,\omega)+S^\eps_a(\lambda\,x,\lambda\,x+(1-\lambda)\,y,\omega).
\]
For the sign of $\phi_a$, we take $v\in\S_c$. From
\eqref{spessore} we know that, for every $a\geq c$, there exists
$\delta_a\geq 0$ with
\[
    S^\eps_a(0,q,\omega)=\eps\, S_a(0,q/\eps,\omega)
    \geq
    \delta_a|q|+\frac{v(q/\eps,\omega)-v(0,\omega)}{1/\eps},
\]
for any $q\in\R^N$. In addition,  $\delta_a>0$ when $a>c$. By the
sublinear character of $v(\cdot,\omega)$ a.s. in $\omega$, we
obtain in the limit $\phi_a(q)\geq \delta_a|q|$. \end{dimo}

%\begin{oss}
%The second part of the above proof  is based on an argument
%provided in \cite{Souga-personal}. It  is also needed in
%\cite{Souga99} to complete the proof of Theorem 1.\\
%\end{oss}

\smallskip

We proceed recalling a result proved in \cite{ReTa00, Souga99}. In
what follows, we  will denote by $L$ the Lagrangian associated
with $H$ via the Fenchel transform, i.e.
\[
L(x,q,\omega):=\max_{p\in\R^N}\big(\,\langle q,p\rangle-
H(x,p,\omega)\big),\qquad
\hbox{$(x,q,\omega)\in\R^N\times\R^N\times\Omega$.}
\]
\begin{prop}\label{prop effective L}
For every $x,y\in\R^N$, $\omega\in\Omega$ and $t>0$
let
\[
h_t(x,y,\omega):=\inf\left\{\int_0^t
L(\gamma,\dot\gamma,\omega)\,\dd
s\,:\,W^{1,1}([0,t],\R^N),\,\gamma(0)=x,\,\gamma(t)=y\,\right\}.
\]
Then there exists a convex and superlinear function $\overline
L:\R^N\to\R$ such that for any $\omega$ in a set $\Omega_0$ of
probability 1 the following convergence holds
\begin{equation}\label{eq cv h_t}
\frac{h_t(0,tq,\omega)}{t}\ \underset{t\to +\infty}{\ucv}\
\overline L(q),\qquad\hbox{$q\in\R^N$.}\medskip
\end{equation}
\end{prop}

The function $\overline L$ is called the {\em effective
Lagrangian}, and the {\em effective Hamiltonian} is accordingly
defined through the Fenchel transform as follows:
\[
\overline H(p)=\max_{q\in\R^N}\left(\,\langle p,q\rangle-\overline
L(q)\right)\qquad\hbox{for every $p\in\R^N$.}
\]

\begin{teorema}\label{teo relation effective H}
For every $a\geq c_f$, the stable norm  $\phi_a$ is the support
function of the $a$--sublevel  of the effective Hamiltonian
$\overline H$.
\end{teorema}

\begin{proof} We denote by $\overline \sigma_a(\cdot)$ the support
function of the $a$--sublevel of $\overline H$. From the
inequality $L(x,q,\omega)+a\geq \sigma_a(x,q,\omega)$, holding for
any $a\geq c_f$ and $(x,q,\omega) \in
\R^N\times\R^N\times\Omega_f$, we infer
\[
\frac{h_t(0,\lambda\,t\,q,\omega)+a\,t}{\lambda\,t}\geq \frac
1{\lambda\,t} S_a(0,\lambda\,t\,q,\omega)
\]
for every $\lambda>0$ and $t>0$. Passing to the limit for $t$
going to $+ \infty$ we find
\[
\lambda^{-1}\,\left(\overline L(\lambda\,q)+a\right)\geq
\phi_a(q)\,
\]
and, taking into account the  identity
\begin{equation*}%\label{eq sigma_a}
\overline\sigma_a(q)=\inf_{\lambda>0}\left\{\lambda^{-1}\,\left(\overline
L(\lambda\,q)+a\right)\right\},
\end{equation*}
we conclude that $\overline\sigma_a(\cdot)\geq \phi_a(\cdot)$. We
divide the proof of the converse inequality in two steps.\\
\indent{\em Case 1: $a>c_f$.} Clearly, it is enough to show that
$\overline\sigma_a(q)\leq \phi_a(q)$ for every $q\in\Sf^{N-1}$.
Let us fix such a $q$ and pick an $\omega_0$ such that both
\eqref{eq cv h_t} and \eqref{eq cv metrics} hold. For every
$n\in\N$, let $\gamma_n:[0,\ell_n]\to\R^N$ be a curve
parameterized by the arc--length with $\gamma_n(0)=0$,
$\gamma_n(\ell_n)=nq$ and such that
\[
1+S_a(0,nq,\omega_0)>\int_0^{\ell_n}
\sigma_a(\gamma_n,\dot\gamma_n,\omega_0)\,\dd t.
\]
We first claim that there exists a constant $C$ such that
$n\leq\ell_n\leq C\,n$ for every $n\in\N$. Indeed, let
$v\in\D{Lip}(\R^N)$ such that $H(x,Dv(x),\omega_0)\leq c_f$ a.e.
in $\R^N$. As $a>c_f$, there exists by \eqref{spessore} a constant
$\delta_a>0$ such that
\[
\int_0^{\ell_n}\sigma_a(\gamma_n,\dot\gamma_n,\omega_0)\,\dd t
    \geq
    v(nq)-v(0)+\delta_a\,\ell_n,
\]
so the claim follows with $C=(1+2\,\kappa_a)/\delta_a$. According
to the results proved in \cite{D1-06} (cf. Lemma 3.4, Proposition
3.7 and Lemma 3.14) there exists a Borel--measurable function
${\underline\lambda\,}_a:\R^N\times\Sf^{N-1}\to [0,+\infty)$ such
that
\[
L(x,{\underline\lambda\,}_a(x,v)\,
{v},\omega_0)=\sigma_a(x,{\underline\lambda\,}_a(x,v)\,
{v},\omega_0)-a\qquad\hbox{for every
$(x,v)\in\R^N\times\Sf^{N-1}$.}
\]
Furthermore ${\underline\lambda\,}_a$ enjoys the following
inequality
\[
\frac{1}{\lambda_a}\leq{\underline\lambda\,}_a(x,{v})\leq
\lambda_a\qquad\hbox{for every $(x,v)\in\R^N\times\Sf^{N-1}$,}
\]
where $\lambda_a$ is a positive real constant depending on $H$ and
$a$ only. Set
\[
{f}_n(s):=\int_0^s
\frac{1}{{\underline\lambda\,}_a(\gamma_n(\varsigma),\dot\gamma_n(\varsigma))}\,\dd
\varsigma\qquad\hbox{for any $s\in[0,\ell_n],$}
\]
and $\varphi_n=f_n^{-1}$ on $[0,f_n(\ell_n)]$. It is easily seen
that $\varphi_n$ is a strictly increasing bi--Lipschitz
homeomorphism from $[0,f_n(\ell_n)]$ to $[0,\ell_n]$, and that
$n/\lambda_a\leq f_n(\ell_n)\leq C\lambda_a\,n$. Arguing as in
\cite{D1-06}, we get that the  curve
\[
\xi_n(s):=\left(\gamma_n\comp\varphi_n\right)(s),\qquad
s\in[0,f_n(\ell_n)]
\]
is a reparameterization of $\gamma_n$ such that
\[
\int_0^{\ell_n}\sigma_a(\gamma_n,\dot\gamma_n,\omega_0)\,\dd t =
\int_0^{f_n(\ell_n)}
\left(L(\xi_n,\dot\xi_n,\omega_0)+a\right)\,\dd t.
\]
For each $n\in\N$, let $f_n(\ell_n)=\lambda_n\,n$ {with
$\lambda_n\in [1/\lambda_a,C\lambda_a]$}. Up to subsequence, we
can assume that $\lambda_n\to\lambda$ as $n\to +\infty$. Then we
get
\[
    \frac{1+S_a(0,nq,\omega_0)}{n}
    \geq
    \frac{h_{\lambda_n n}(0,nq,\omega_0)}{n}+\lambda_n\,a
    =
    \lambda_n\,
    \frac{h_{\lambda_n n}(0,\lambda_n n\,q/\lambda_n,\omega_0)}{\lambda_n n}+\lambda_n\,a
\]
and sending $n\to +\infty$ we finally obtain
\[
\phi_a(q)\geq \lambda\left(\overline L(q/\lambda)+a\right)\geq
\overline\sigma_a(q).
\]
\smallskip\\
\indent{\em Step 2: $a=c_f$.} We want to show that
$\phi_{c_f}(\cdot)\equiv \overline\sigma_{c_f}(\cdot)$. By the
previous step and by definition of $\overline\sigma_a$ we have
\[
\overline\sigma_{c_f}(q)=\inf_{a>{c_f}}
\overline\sigma_a(q)=\inf_{a>{c_f}} \phi_a(q)\qquad\hbox{for every
$q\in\R^N$.}
\]
We therefore get the assertion showing
\[
\phi_{c_f}(q)=\inf_{a>{c_f}} \phi_a(q) \quad\text{ for every
$q\in\R^N$.}\]
 The inequality $\phi_{c_f}(q)\leq\inf_{a>{c_f}}
\phi_a(q)$ comes directly  from the monotonicity   of $a\mapsto
S_a(0,nq,\omega)$. For  the converse, we fix $q$ and invoke
Kingman's Subadditive Ergodic Theorem to get
\begin{equation*}%\label{eq subadditive ergodic thm}
\phi_a(q)=\inf_{n\in\N}\, \frac 1n \, \EE\big(
S_a(0,nq,\cdot)\big)\qquad\hbox{for every $a>c_f$}.
\end{equation*}
Since $S_{c_f}(0,nq,\omega)=\inf_{a>{c_f}} S_a(0,nq,\omega)$ for
every $\omega$, the Monotone Convergence Theorem then implies
\[ \frac 1n \,
    \EE( S_{c_f}(0,nq,\cdot)
    =
    \frac 1n \, \inf_{a>{c_f}} \EE(S_{a}(0,nq,\cdot))
    \geq
    \inf_{a>{c_f}} \phi_a(q),
\]
and sending $n\to +\infty$ we obtain
$\phi_{c_f}(q)\geq\inf_{a>{c_f}} \phi_a(q)$, as
claimed.\end{proof} \medskip
\end{section}

\begin{section}{Final  results} \label{main}

In this section we exploit the previous analysis  on stable norms
and the Lax formula given in Section \ref{sez HJ}  to establish
relevant properties of the effective Hamiltonian as well
as some existence/nonexistence result for exavt and approximate correctors. \smallskip \par

\begin{teorema}\label{teo min effective H}
$\min_{\R^N} \overline H=c_f$.
\end{teorema}

\begin{dimo}
The inequality $\min_{\R^N} \overline H\leq c_f$ is immediate
since the $c_f$--sublevel of $\overline H$ is nonempty, being
$\overline\sigma_{c_f}= \phi_{c_f}$ finite by Theorems \ref{teo
metric homogenization} and \ref{teo relation effective H}.

Let us prove $\min_{\R^N} \overline H\geq c_f$. Pick an $\omega$
in $\Omega_f$. By Proposition 2-1.1 in \cite{CoIt00} we know that,
for every $a<c_f$,
\[
\inf_{t>0}\,\left( h_t(0,0,\omega)+at\right)=-\infty,
\]
in particular there exists $t_0>0$ such that
$h_{t_0}(0,0,\omega)+at_0<0$. Since $h_{t+s}(0,0,\omega)\leq
h_t(0,0,\omega)+h_s(0,0,\omega)$ for every $s,t>0$ by the
definition of $h_t$, we infer
\[
    \liminf_{t\to +\infty}\,\frac{h_t(0,0,\omega)+at}{t}
    \leq
    \lim_{n\to +\infty}\,\frac{n \left(h_{t_0}(0,0,\omega)+a t_0\right)}{n t_0}
    <
    0.
\]
In view of Proposition \ref{prop effective L}, we get $\overline
L(0)+a<0$ for every $a<c_f$, that is $-\overline L(0)\geq c_f$.
The assertion follows since $-\overline L(0)=\min_{\R^N} \overline
H$.
\end{dimo}

\medskip

\ \\
\indent  We exploit  Proposition \ref{prop Lax in S} to get

\begin{prop}\label{prop stationary subsolutions}
Let $a\geq c_f$ such that the corresponding stable norm is
nondegenerate. Then equation \eqref{eq HJa} admits admissible
stationary subsolutions.
\end{prop}

\begin{dimo}
By hypothesis  there exists $\delta_a>0$ such that $\phi_a(q)\geq
\delta_a|q|$ for every $q\in\R^N$. By Theorems \ref{teo metric
homogenization} and \ref{teo relation effective H} we derive
\[
\liminf_{|y|\to +\infty}\frac{S(y,x,\omega)}{|y-x|}\geq
\delta_a\qquad\hbox{for every $x\in\R^N$}
\]
a.s. in $\omega$. In particular $\inf_{y\in\R^N} S_a(y,x,\omega)>
- \infty$ a.s. in $\omega$ for every fixed $x\in\R^N$. According
to Proposition \ref{prop Lax in S}, an admissible stationary
subsolution of \eqref{eq HJa} is obtained via \eqref{def u} with
$g(\cdot,\omega)\equiv 0$ and $C(\omega)=\R^N$ for every $\omega$.
\end{dimo}

\ \medskip\\
\indent We make use of the previous results to give a
characterization of the stationary critical value.

\begin{teorema}\label{teo characterize c}
$c=\inf\{a\geq c_f\,:\, \phi_a \;\text{is nondegenerate}\,\}$. If
$c>c_f$ then $\phi_c$ is degenerate but nonnegative.
\end{teorema}

\begin{dimo}
Let us call $\mu$ the infimum appearing in the statement.
According to Theorem \ref{teo metric homogenization},
 $c\geq \mu$. The
converse inequality is apparent by Proposition \ref{prop
stationary subsolutions} since equation \eqref{eq HJa} admits
admissible subsolutions for every $a>\mu$.\par

Assume by contradiction that $c >c_f$ and $\phi_c$ is
nondegenerate. Then, since $\phi_a$ coincide with $a$--sublevel of
$\overline H$, by continuity the same property holds for $\phi_a$
with $a<c$ and suitably close to $c$. That is in contradiction
with Proposition \ref{prop stationary subsolutions}.
\end{dimo}

\begin{oss}\label{oss characterize c}
In view of Theorem \ref{teo relation effective H}, the property of
being $\phi_a$ non--degenerate is equivalent to the condition
$0\in\D{Int}(\overline Z_a)$. So Theorem \ref{teo characterize c}
can be equivalently restated by saying that $c=\inf\{a\geq
c_f\,:\,0\in\D{Int}(\overline Z_a)\,\}$, with $0\in\partial
\overline Z_c$ whenever $c>c_f$. See \cite{DavSic08} for the
analogy with the 1--dimensional case.
\end{oss}
\indent We derive from  Proposition \ref{prop stationary
subsolutions} and Theorem \ref{teo characterize c}
\begin{equation*}\label{def c1}
c:=\inf\{a\in\R\,:\,\text{\abra{eq HJa} admits stationary
subsolutions} \}.
\end{equation*}
The infimum appearing in the formula is not necessarily a minimum,
namely we cannot expect, in general, to find stationary critical
subsolutions.\medskip

The next theorem relates the effective Hamiltonian to the
stationary critical value.  The result has been already
established by Lions and Souganidis in \cite{LiSou03} through PDE
techniques. We propose a new, simpler proof based on the
properties of the intrinsic metrics and on Theorems \ref{teo
metric homogenization}, \ref{teo relation effective H}.

\begin{teorema} \label{effettivo} $ \overline H$  coincides with the function associating
to any $P \in \R^N$ the stationary critical value of the
Hamiltonian $H(x, P+p, \omega)$.
\end{teorema}

\begin{dimo}
We first prove the assertion for $P=0$, i.e., with the notation
used so far, that $\overline H(0)=c$. We know by Theorem \ref{teo
metric homogenization}  that $\phi_c(q)\geq 0$ for every
$q\in\R^N$. This means, in view of Theorem \ref{teo relation
effective H}, that $0\in\overline Z_c$, i.e. $c\geq \overline
H(0)$. On the other hand the inequality $c>\overline H(0)$, i.e.
$0\in\D{Int}(\overline Z_c)$, may occur only  when $c=c_f$ by
Remark \ref{oss characterize c}, but this is not possible since
$\overline H(0)\geq c_f$ by Theorem \ref{teo min effective H}.

To prove the assertion for any $P\in\R^N$, we apply the previous
argument to the Hamiltonian $H_P(x,p,\omega):=H(x,P+p,\omega)$ and
derive that $\overline H_P(0)=c_P$, where $\overline H_P$ and
$c_P$  are the associated effective Hamiltonian and stationary
critical value. To get the assertion, it is left to show that
$\overline H_P(0)=\overline H(P)$. To this purpose, denote by
$L_P$ the Lagrangian associated with $H_P$. It is easily seen that
$L_P(x,q,\omega)$ coincides with $L(x,q,\omega)-\langle
P,q\rangle$, so the claimed equality follows from the definition
of effective Hamiltonian by exploiting  Proposition \ref{prop
effective L} with $L_P$ in place of $L$.
\end{dimo}

\ \smallskip\\
\indent We now address our attention to the issue of the
existence/nonexistence  of exact  or approximate correctors.

\begin{teorema}\label{teo correctors}\ \vspace{1ex}
\begin{itemize}
    \item[{\em (i)}] If $c=c_f$ and  $\A_f(\omega)\not =\emptyset$ a.s. in $\omega$,
    then there exists a corrector
    for \eqref{eq critica}.\smallskip
    \item[{\em (ii)}] If $\phi_c$ is nondegenerate, then
     $c=c_f$ and a corrector for \eqref{eq critica} exists if and
    only if $\A_f(\omega)\not=\emptyset$ a.s. in $\omega$. \smallskip
    \item[{\em (iii)}] If $\phi_c$ is nondegenerate and  $\A_f(\omega) \neq \emptyset$ a.s. in
    $\omega$,
    then the classical Aubry set  is an uniqueness  set for
    the critical equation, in the sense that two correctors
    agreing on $\A_f(\omega)$ a.s. in $\omega$, coincide on the
    whole $\R^N$ a.s. in $\omega$. More precisely, any corrector
    $u$ can be written as
    \begin{equation}\label{laxancora}
    u(x,\omega)= \inf \{u(y,\omega)+S_c(y,x,\omega) \,:\,y \in
    \A_f(\omega)\} \quad\text{a.s. in $\omega$}.
\end{equation}
\end{itemize}
\end{teorema}

\begin{dimo}
{\em (i)} \ By Proposition \ref{randomclosed} $\A_f(\omega)$ is a
stationary closed random set. Hence, for any $g\in\S_c$, the
function $u$ given by the Lax formula with $a=c$, $\A_f(\omega)$
as source set  and  trace $g$ on it, is a corrector by   Theorem \ref{teo 6.7} (i).\smallskip\\
\indent{\em (ii)}\quad  The equality $c=c_f$ follows from Theorem
\ref{teo characterize c}. If $\A_f(\omega)\not=\emptyset$ a.s. in
$\omega$, a corrector for \eqref{eq critica} exists by assertion
(i). To prove the converse implication, let us assume by
contradiction that a corrector $u$ does exist and that
$\A_f(\omega)=\emptyset$ a.s. in $\omega$. Pick $\omega$ such that
$\A_f(\omega)=\emptyset$ and $u(\cdot,\omega)$ is a viscosity
solution of \eqref{eq critica}. Take $n\in\N$. Since the classical
Aubry set is empty, we derive by Theorem \ref{teo 6.7} that
$u(\cdot,\omega)$ is the unique viscosity solution of the
 Dirichlet Problem
\begin{eqnarray*}
\begin{cases}
H(x,D\phi(x),\omega)=c_f&\qquad\hbox{in $B_n$}\\
\phi(x)=u(x,\omega)&\qquad\hbox{on $\partial B_n$,}
\end{cases}
\end{eqnarray*}
and
\[
u(x,\omega)=\min_{y\in\partial
B_n}\{u(y,\omega)+S_{c_f}(y,x,\omega)\,\}\qquad\hbox{$x\in B_n$}.
\]
We deduce that there exists a {diverging} sequence $(y_n)_n$ such
that
\[
S_{c_f}(y_n,0,\omega)=u(0,\omega)-u(y_n,\omega)\qquad\hbox{for
every $n\in\N$.}
\]
By exploiting the fact that $u(\cdot,\omega)$ is sublinear a.s. in
$\omega$, we derive
\[
\min_{q\in\Sf^{N-1}}\phi_{c_f}(q)=\liminf_{|y|\to
+\infty}\frac{S_{c_f}(y,0,\omega)}{|y|} \leq \liminf_{n\to
+\infty}\frac{u(0,\omega)-u(y_n,\omega)}{|y_n|}=0\qquad\hbox{a.s.
in $\omega$,}
\]
in contrast with the hypothesis that $\phi_{c_f}$ is
nondegenerate. \smallskip\\
\indent{\em (iii)}\quad We  take a corrector $u$ and fix an
$\omega \in \Omega_f$ such that $u(\cdot,\omega)$ is a solution to
\eqref{eq critica} sublinear at infinity  and $\A_f(\omega) \neq
\emptyset$. Arguing as for item (ii), we see that for any given
$x_0 \in \R^N$ there is  $n=n(\omega)$ such that
\[ u(x_0,\omega) < \inf\{u(z,\omega) + S_c(z,x_0,\omega)\,:\, z - x_0 \in
\partial B_n\}.\]
 According to Theorem \ref{teo 6.7} (iii), we deduce the
existence of $y_0 \in \A_f(\omega)$ with
\[u(x_0,\omega) =u(y_0,\omega)  + S_c(y_0,x_0,\omega),\]
and consequently
\[ u(x_0,\omega) \geq \inf\{u(y,\omega) + S_c(y,x_0,\omega)\,:\, y  \in
\A_f(\omega)\}.\] On the other side by Theorem \ref{teo 6.7} (i)
the right hand--side of the previous formula is the maximal
subsolution to \eqref{eq critica} taking the value
$u(\cdot,\omega)$ on $\A_f(\omega)$, which  implies
\eqref{laxancora}.
\end{dimo}

\ \smallskip\\
\indent In the case where $c=c_f=\sup_{x \in \R^N} \min_{p \in \R^N}
H(x,p,\omega)$ a.s. in $\omega$, item (i) of the previous theorem
can be complemented as follows:

\begin{prop}\label{equi} Assume  $c=c_f=\sup_{x \in \R^N} \min_{p \in \R^N}
H(x,p,\omega)$  and $\A_f(\omega) = \emptyset$ a.s. in $\omega$,
then  equation \eqref{eq critica} admits approximate correctors.
\end{prop}

\begin{dimo}
 We fix $g \in \S_c$ and  define, for any  $\delta>0$
\[\Eq_\delta(\omega)=\{x\in\R^N\,: \min_{p \in \R^N}
H(x,p,\omega) \geq c - \delta\}.\] Arguing as we did for $\Eq$ in
Proposition \ref{randomclosed}, we see that $\Eq_\delta(\omega)$
is a closed random stationary set and is, in addition, a.s.
nonempty. We claim that Lax formula with $\Eq_\delta(\omega)$ as
source set and $g$ as trace on it, provides  a
$\delta$--approximate corrector.\par

We denote by $v_\delta$ the random function constructed as above
indicated.  By Proposition \ref{prop Lax in S} we already know
that $v_\delta$ is an admissible subsolution to \eqref{eq
critica}, and a solution as well on
$\R^N\setminus\Eq_\delta(\omega)$ a.s. in $\omega$. Further, if
$\phi$ is a $C^1$ test function touching $v_\delta(\cdot,\omega)$
at $y \in \Eq_\delta(\omega)$  from below, then
$H(y,D\phi(y,\omega),\omega)\geq c -\delta$ by the very definition
of $\Eq_\delta(\omega)$, as it was to be shown.\end{dimo}

\ \smallskip\\
\indent A class of critical stochastic equations  satisfying  the
assumptions of the previous theorem are those of Eikonal type
\begin{equation}\label{eq Eikonal}
|Du(x,\omega)|^2=V(x,\omega)^2\qquad\hbox{in $\R^N$,}
\end{equation}
where $V:\R^N\times\Omega\to\R$ is a jointly measurable function
satisfying:\smallskip
\begin{itemize}
    \item[(1)]\quad $V(x+z,\omega)=V(x,\tau_z\omega)$ \quad for every
$x,z\in\R^N$ and $\omega\in\Omega$;\medskip
    \item[(2)]\quad
$V(\cdot,\omega)$ \ is continuous on $\R^N$\quad for every
$\omega$;\medskip
    \item[(3)]\quad $0=\inf_{\R^N} V(\cdot,\omega)<\sup_{\R^N}
V(\cdot,\omega)<+\infty$ \quad a.s. in $\omega$;\medskip
    \item[(4)]\quad $V(x,\omega)>0$ for every $x\in\R^N$ \quad a.s. in
$\omega$.\medskip
\end{itemize}

We can show, by exploiting an example in \cite{LiSou03}, that a
random function of this type does always exist in any space
dimension, cf. Example \ref{esempio degenere} below.

\medskip

\indent  If we add to the assumptions of Proposition \ref{equi}
the nondegeneracy  of the critical stable norm, then we can also
assert, according to Theorems \ref{teo correctors} (ii),  the
nonexistence  of exact correctors. In dimension 1 a sufficient
conditions for such a nondegeneracy  is that there is a strict
critical admissible subsolution; this role is for instance played
by the null function for the above described Eikonal class.  The
point is that in the 1--dimensional case  the $a$--sublevel of the
effective Hamiltonian, for any $a \geq c_f$, coincides with the
averaged $a$--sublevel  $\EE[Z_a]$ of $H$,  given by
\[ \EE[Z_a]= \{\EE(\Phi) \,:\, \Phi \;\text{measurable selection of
$\omega \mapsto Z(0,\omega)$} \}.\]
 Due to the stationarity properties of the Hamiltonian, the
 previous definition
 does not change if we replace $0$ by any other $x \in \R^N$. For
 every $x$, the map $\omega \mapsto Z(x,\omega)$ is a random closed
 set taking compact convex values and (see \cite{Molchanov})
 \[\sigma_{\EE[Z_a]}(q)= \EE\big(\sigma_a(x,q,\omega)\big),
 \qquad  q \in \R^N,\]
 where $\sigma$ indicates the support function.
 In dimension $1$  any such selection with vanishing mean gives rise, by integration, to
an admissible subsolution, see \cite{DavSic08}, which  actually  explains why the existence of a strict
admissible critical subsolution implies the nondegeneracy of the critical
stable distance.

The situation is quite different in the multidimensional setting,
where we can just assert that the  sublevels of $\overline H$
 are contained in the corresponding averaged sublevels of $H$. The next
example shows that this inclusion can be strict and that the
critical stable norm can be (even completely) degenerate  in
presence of a strict admissible critical subsolution.

\begin{esempio}\label{esempio degenere}
We provide below an example in dimension $N=2$ of a function $V$
satisfying assumptions (1)--(4) such that, for every $\omega$ in a
set of probability 1,
\[
\lim_{t\to +\infty}\frac{S_0(0,tq,\omega)}{t}=0\qquad\hbox{for
every $q\in\Sf^{N-1}$,}
\]
where $S_0$ is the distance associated with
$H(x,p,\omega):=|p|^2-V(x,\omega)^2$ via \eqref{eq S} with $a=0$.
According to the results obtained in the previous section, we
derive that the  corresponding stable norm is null, i.e.
completely degenerate. Note that the null function is a strict
admissible critical subsolution. \par
 To this purpose, we start by defining a function
$V_0:\T^2\to\R$ as follows:
\[
V_0(z_1,z_2)=2-\cos(2\pi z_1)-\cos(2\pi z_2)\qquad
(z_1,z_2)\in\T^2.
\]
Let us choose a $\lambda\in\R\setminus\Q$ and set $v_1=(1,0),\
v_2=(\lambda,0),\ v_3=(0,1),\ v_4=(0,\lambda).$ Next we choose as
$\Omega$ the torus $\T^4$, as $\PP$ the Lebesgue measure
restricted to $\T^4$, and as $\F$ the $\sigma$--algebra of Borel
subsets of $\Omega$. We define a group $(\tau_x)_{x\in\R^2}$ of
translations on $\Omega$ as follows:
\[
\big(\tau_x(\omega)\big)_i\equiv\omega_i+\langle
v_i,x\rangle\quad\hbox{(mod 1)}\qquad\hbox{for every
$\omega\in\Omega$ and $i\in\{1,2,3,4\}$.}
\]
The group of translations $(\tau_x)_{x\in\R^2}$ is ergodic, see
for instance Appendix A in \cite{DS1-08}. We define a function $V$
on $\R^2\times\Omega$ as
\[
V(x,\omega)=V_0(\omega_1+x_1,\omega_2+\lambda
x_1)\,V_0(\omega_3+x_2,\omega_4+\lambda x_2),
\]
where {$x=(x_1,x_2)\in\R^2$ and
$\omega=(\omega_1,\omega_2,\omega_3,\omega_4)\in\Omega$.} Clearly,
$V$ is a jointly measurable function satisfying the above
assumptions (1)--(3). Furthermore, $V$ verifies assumption (4).
Indeed, the function $V(\cdot,\omega)$ attains its infimum on
$\R^2$ (that is, the value 0) if and only if
$\omega\in\big(\Sigma\times\T^2\big)\cup\big(\T^2\times\Sigma\big)$,
where
\[
\Sigma:=\{(z_1,z_2)\in\T^2\,:\,V_0(z_1+t,z_2+\lambda t)=0\
\hbox{for some $t\in\R$\ }\}.
\]
We claim that $|\Sigma|=0$. More precisely,
\[
\Sigma=\{(z_1,z_2)\in\T^2\,:\,\lambda z_1-z_2=\lambda n-m\
\hbox{for some $n,m\in\Z$\,}\}.
\]
To see this, note that $(z_1,z_2)\in\Sigma$ if and only if there
exists $t\in\R$ such that
\[
z_1+t=n,\ z_2+\lambda t=m\quad\hbox{for some $n,m\in\Z$,}
\]
and this occurs if and only if $\lambda z_1-z_2=\lambda n-m$ for
some $n,m\in\Z$, as it was claimed.

According to the results of Section \ref{stable}, we know that,
for every $\omega$ in a set $\Omega_0$ of probability 1, we have
\begin{equation}\label{eq hom}
\overline\sigma_0(q)=\lim_{t\to
+\infty}\frac{S_0(x,x+tq,\omega)}{t}\qquad\hbox{for every
$x\in\R^N$ and $q\in\Sf^{N-1}$.}
\end{equation}
To prove that $\overline\sigma_0(\cdot)\equiv 0$, it is enough, by
the properties enjoyed by $\overline\sigma_0$, to show that
$\overline\sigma_0(\mathbf{e_1})\leq 0$ and
$\overline\sigma_0(\mathbf{e_2})\leq 0$, where
$\{\mathbf{e_1},\mathbf{e_2}\}$ is the canonical basis of $\R^2$.
We only prove the assertion for $\mathbf{e_1}$, being the other
case analogous. Pick an
$\omega=(\omega_1,\omega_2,\omega_3,\omega_4)\in\Omega_0$, and fix
an $\eps>0$. The orbit
$\{(\omega_3+x_2,\omega_4+\lambda\,x_2)\,:\,x_2\in\R\,\}$ is dense
in $\T^2$ for $\lambda$ is irrational, so there exists a point
$x_\eps=(0,y_\eps)\in\R^2$ such that
$V_0((\omega_3,\omega_4),(0,y_\eps))<\eps$. By  moving along the
segment joining $x_\eps$ to $x_\eps+t\mathbf{e_1}$, we obtain
\[
    S_0(x_\eps,x_\eps+t\mathbf{e_1},\omega)
    \leq
    \int_0^1 V(x_\eps+s\,t\mathbf{e_1},\omega)|t\,\mathbf{e_1}|\,\dd s
    \leq
    4\eps\,t,
\]
from which we derive $\overline\sigma_0(\mathbf{e_1})\leq 4\eps$
by \eqref{eq hom}. This concludes the proof as $\eps$ is an
arbitrarily chosen positive number.\\
\end{esempio}
\end{section}

\begin{appendix}
\begin{section}{}
\indent The aim of this section is to provide complete proofs of
Theorems \ref{ascoli}, \ref{teorema mean 0} and \ref{teorema
constant mean}. For this, we need to develop some preliminary
material.\smallskip\par

\indent  We define an $N$--parameter group $(U_x)_{x\in\R^N}$ of
isometries on $\left(L^2(\Omega)\right)^N$
 via
\begin{eqnarray*}
U_x:\left(L^2(\Omega)\right)^N&\to& \left(L^2(\Omega)\right)^N\\
\Phi(\omega)\qquad&\mapsto&\Phi(\tau_x\omega)
\end{eqnarray*}
for every $x\in\R^N$. The group $(U_x)_{x\in\R^N}$ is strongly
continuous, in the sense that
\[
\lim_{x\to
0}\|U_x\Phi-\Phi\|_{\left(L^2(\Omega)\right)^N}=0,\qquad{\Phi\in
\left(L^2(\Omega)\right)^N,}
\]
see \cite{JiKoOl}. Using this property, it is easy to prove:

\begin{lemma}
Let $\gamma:[0,T]\to\R^N$ be of class $C^1$. Then for every
$\Phi\in \left(L^2(\Omega)\right)^N$ the map $t\mapsto \langle
U_{\gamma(t)}\Phi,\dot\gamma(t)\rangle$, from $[0,T]$ to
$L^2(\Omega)$, is continuous.
\end{lemma}

Hence, for every curve $\gamma:[0,T]\to\R^N$ of class $C^1$ we can
give a meaning, in the Cauchy sense,  to the integral
\begin{equation}\label{def int}
\int_0^T \langle U_{\gamma(t)}\Phi,\dot\gamma(t)\rangle\,\dd t =
\int_0^T  \langle \Phi(\tau_{\gamma(t)} \omega) ,
\dot\gamma(t)\rangle\,\dd t
\end{equation}
as an element of $L^2(\Omega)$. We note that \eqref{def int} is
invariant under changes of parameterization; moreover, it makes
sense even when $\gamma$ is piecewise $C^1$, i.e. it is continuous
and of class $C^1$ on $[0,T]$ up to a finite  set of points. For
any such $\gamma$ we write
\[
\int_\gamma \Phi(\omega) :=\int_0^T \langle
U_{\gamma(t)}\Phi(\omega),\dot\gamma(t)\rangle\,\dd
t\qquad\hbox{a.s. in $\omega$.}
\]
The following result holds:

\begin{lemma}\label{lemma continuita}
Let $\gamma:[0,T]\to\R^N$ be a piecewise $C^1$ curve. Then the map
\[
\Phi\mapsto \int_\gamma \Phi
\]
is linear and continuous from $\left(L^2(\Omega)\right)^N$ to
$L^2(\Omega)$.
\end{lemma}
It is a direct consequence of the previous lemma that, if $\Phi_n
\rightarrow \Phi$ in $(L^2(\Omega))^N$, then,  for any given
piecewise $C^1$ curve $\gamma$, one can extract a subsequence
$\Phi_{n_k}$ with
\[ \int_\gamma \Phi_{n_k}(\omega) \rightarrow \int_\gamma \Phi(\omega)
\quad\text{a.s. in $\omega$}.
\]
Note that here the subsequence depends on the curve $\gamma$.  The
step forward in the next result is to show that, under suitable
additional assumptions, the sequence $(n_k)_k$ can be chosen in
such a way that the above convergence takes place for any curve.

\begin{lemma}\label{rotto} Let $(\Phi_n)_n$ be a sequence in
$(L^\infty(\Omega))^N$ with  $\sup_n\|\Phi_n\|_\infty  \leq
\kappa$ for some $\kappa >0$. If $\Phi_n$ converges in
$(L^2(\Omega))^N$ to some function $\Phi$, then, up to extraction
of a subsequence,
\[
\int_\gamma \Phi_n(\omega) \rightarrow \int_{\gamma} \Phi(\omega)
\quad\text{a.s. in $\omega$}
\]
for all piecewise $C^1$ curve $\gamma$.
\end{lemma}
\begin{dimo}  Up to extraction of a subsequence,   $\Phi_n$ a.s. converges to $\Phi$. The set
\[
\{(x,\omega)\,:\, \lim_n \Phi_n(\tau_x \omega) =\Phi(\tau_x
\omega)\,,\, |\Phi_n(\tau_x \omega)| \leq \kappa \;\text{for any
$n$}\},
\]
is clearly measurable with respect to the product
$\sigma$--algebra $\mathcal B(\R^n) \otimes \F$, and so  the
almost sure convergence of $\Phi_n$ to $\Phi$  and the boundedness
assumption on $\Phi_n$ imply that its $x$--sections  have
probability 1 for any fixed $x\in\R^N$. We derive from Fubini's
Theorem  that there exists a set $\Omega'$ of probability $1$ such
that, if $\omega\in\Omega'$
\begin{equation}\label{ascoli1}
    \Phi_n(\tau_x \omega) \rightarrow \Phi(\tau_x \omega),\quad
    \sup_n|\Phi_n(\tau_x \omega)| \leq \kappa \qquad\text{for a.e. $x\in\R^N$.}
\end{equation}
 The set
\[X(\omega) = \{x \, :\, \Phi_n(\tau_x \omega) \rightarrow \Phi(\tau_x \omega),
 \quad\sup_n|\Phi_n(\tau_x \omega)| \leq \kappa\}
\]  is accordingly a stationary random set with volume fraction
$1$. Therefore, given a piecewise $C^1$ curve $\gamma$, by
applying the Robin's Theorem  with $\mu$ equal to $\hh^1$
restricted to $\gamma$, we see that $\hh^1\left(\gamma\setminus
X(\omega)\right)=0$ for $\omega$ belonging to some subset of
$\Omega$ of probability $1$. For such an $\omega$  the claimed
convergence on $\gamma$  holds true thanks to the Dominated
Convergence Theorem.
\end{dimo}

\medskip

We proceed by giving  a closer look to   the differentiability
properties of Lipschitz random function with stationary
increments.

\begin{prop}\label{dif} Let $v$ be a Lipschitz random function
with stationary increments. Then there exists $\Phi \in
(L^\infty(\Omega))^N$  such that:
\begin{itemize}
    \item[{\em (i)}] for every $\omega$ in a set of probability 1
    \begin{equation}\label{dif01}
 D v(x,\omega) = \Phi(\tau_x \omega) \quad\text{for a.e. $x \in \R^N$;}
    \end{equation}
    \item[{\em (ii)}] for every closed piecewise $C^1$ curve
    $\gamma$
    \[
    \int_\gamma \Phi(\omega) = 0\qquad\hbox{ a.s. in
$\omega$}.
    \]
\end{itemize}
In addition, the  equality \abra{dif01}  holds, for {\em any}
fixed $x$, a.s. in $\omega$.
\end{prop}

\begin{dimo}
Let $\kappa$ be a positive constant such that $v(\cdot,\omega)$ is
$\kappa$--Lipschitz for every $\omega$. For each $i\in\{1, \cdots,
N\}$ we define
\[
w_i(x,\omega) =\sup_{n \in \N}\, \inf_ {h\in \Q \cap B_{1/n}}
\frac{v(x +h\,e_i,\omega)- v(x,\omega)}h.
\]
Such function is jointly measurable in $\R^N \times \Omega$,
satisfies $|w_i(x,\omega)|\leq\kappa$ for every $x\in\R^N$ and
$\omega\in\Omega$, and, in addition, it is stationary, being $v$
with stationary increments. By \abra{staz}, there exists a set
$\Omega'$ of probability $1$ such that
\begin{equation*}
    w_i(x,\omega)= w_i(0,\tau_x \omega) \quad\text{for
a.e. $x\in\R^N$ and any $\omega \in \Omega'$.}
\end{equation*}
Moreover, for every $\omega\in\Omega$,
\begin{equation*}
    w_i(x,\omega)= \partial_{  x_i} v(x,\omega)
  \quad \hbox{for every $x\in\Delta_v(\omega)$.}
\end{equation*}
Since $\Delta_v(\omega)$ has full measure in $\R^N$ by the
Lipschitz character of $v(\cdot,\omega)$, the equality
\abra{dif01} is obtained by setting
$\Phi_i(\cdot)=w_i(0,\cdot)$.\par

For any $\omega \in \Omega$ we set
\[
X(\omega)=\{ x \in \Delta_v(\omega)\,:\,
Dv(x,\omega)=\Phi(\tau_x\omega)\}.
\]
By taking into account that $\Delta_v$ is a stationary random set
(cf. Proposition \ref{prop Delta_v random}) and that $v$ has
stationary increments, we see that $X$ is a random stationary set.
Furthermore, from \abra{dif01} we deduce that its volume fraction
is $1$. By Fubini's Theorem, we deduce that for any fixed
$x\in\R^N$ the equality \eqref{dif01} holds a.s. in $\omega$.

Given a  piecewise $C^1$ closed curve $\gamma:[0,T]\to\R^N$,  we
invoke Robbins' Theorem with $\mu$ equal to $\hh^1$ restricted to
$\gamma$ to deduce that $\hh^1\left(\gamma\setminus
X(\omega)\right)= 0$ for $\omega$ belonging to some subset of
$\Omega$ with probability $1$.  For such an $\omega$ we get
\[
    \int_\gamma \Phi(\omega)
    =
    \int_0^T \langle
Dv({\gamma(t)},\omega),\dot\gamma(t)\rangle\,\dd t
    =
    v(\gamma(T),\omega)-v(\gamma(0),\omega)
\]
and the assertion follows as $\gamma(T)=\gamma(0)$.\medskip
\end{dimo}

Conversely, we have

\begin{prop}\label{teo null rot} Let $\Phi\in
\left(L^\infty(\Omega)\right)^N$ with
\[
   \int_\gamma \Phi(\omega) =0\qquad\hbox{ a.s. in  $\omega$}
\]
for every closed piecewise $C^1$  curve. Then  there exist  a
Lipschitz random  function $v$ with stationary increments  and a
set $\Omega_0$ of probability $1$ such that, for any
 $\omega\in \Omega_0$,
\[
Dv(x,\omega)=\Phi(\tau_x\omega)\quad\hbox{ for a.e.
$x\in\R^N$.\smallskip}
\]
\end{prop}
 \begin{dimo}   Let us set
\[
v(x,\omega)=\int_0^1 \langle \Phi(\tau_{tx}\,\omega),x\rangle\,\dd
t,\qquad\hbox{$(x,\omega)\in\R^N\times\Omega.$}
\]
Note that  $v$
 is jointly measurable in $(x,\omega)$. By  assumption we also
 have that, for every fixed $x\in\R^N$,
\begin{equation}\label{rot1}
   v(x,\omega)=\int_0^T \langle
\Phi(\tau_{\gamma(t)}\,\omega),\dot\gamma(t) \rangle\,\dd
t,\qquad\text{a.s. in $\omega$}
\end{equation}
whenever $\gamma$ is a piecewise $C^1$ curve,  defined in some
interval $[0,T]$, joining $0$ to $x$.   For every $i = 1, \cdots,
N$ we
 derive from \abra{rot1}
 \[
\frac{v(x+h\,e_i,\omega)-v(x,\omega)}{h}=\int_0^1
\Phi_i(\tau_{x+t\,h\,\e_i}\,\omega)\,\dd t\qquad\hbox{for any $x$,
a.s. in $\omega$,}
\]
where $h$ is a  discrete positive parameter. Hence
\begin{eqnarray*}
     \int_\Omega \left| \frac{v(x+h\,e_i,\omega)-v(x,\omega)}{h}
    -\Phi_i(\tau_x\omega) \right |^2\,\dd
   \PP
    =
    \int_\Omega \left| \int_0^1
    \Phi_i(\tau_{x+t\,h\,e_i}\,\omega) -\Phi_i(\tau_x\omega)\,\dd
    t\right|^2\,\dd \PP\\
    \leq
    \int_\Omega \int_0^1 \left|\Phi(\tau_{x+t\,h\,e_i}\,\omega) -\Phi(\tau_x\omega)\right|^2
    \,\dd t\,\dd \PP
    =
    \int_0^1
    \|U_{t\,h\,e_i}\Phi-\Phi\|^2_{\left(L^2(\Omega)\right)^N}\,\dd t,
\end{eqnarray*}
where in the last equality we have used Fubini's Theorem and the
fact that the probability measure $\PP$ is invariant under
$\tau_x$. The previous relation implies, by the strong continuity
of the group $(U_y)_{y\in\R^N}$
\[
\frac{v(x+h\,e_i,\omega)-v(x,\omega)}{h}\underset{h\to 0}\tto
\Phi_i(\tau_x\omega)\qquad\hbox{in $L^2(\Omega)$,}
\]
for every $x\in\R^N$; accordingly, thanks to Fubini's Theorem
\[
\frac{v(x+h\,e_i,\omega)-v(x,\omega)}{h}\underset{h\to 0}\tto
\Phi_i(\tau_x\omega)\qquad\hbox{in $L^2_{loc}(\R^N)$,}
\]
a.s. in $\omega$. This is, in turn, equivalent to the equality
$\partial_{x_i} v(x,\omega)=\Phi_i(\tau_x\omega)$, for $i=1,
\cdots, N$,  in the sense of distributions, a.s. in $\omega$.
Therefore, being $\Phi$ essentially bounded, $v(\cdot,\omega)$ is
Lipschitz--continuous with $Dv(x,\omega)= \Phi(\tau_x\omega)$ for
a.e. $x\in\R^N$ and every $\omega$ in a set $\Omega_0$ of
probability 1.  By suitably assigning the value of
$v(\cdot,\omega)$ on $\Omega\setminus\Omega_0$, we extend the
Lipschitz character of such function to all $\omega$. Therefore
$v$, being also jointly measurable in $\R^N \times \Omega$, is a
Lipschitz random function, as asserted.\par

It is left  to show that it has stationary increments. Let us fix
$z\in\R^N$. Then, by \abra{rot1}, for every $x,y\in\R^N$
\begin{eqnarray}\label{eq incr staz}
v(x+z,\omega)&\!\!\!\!\!-&\!\!\!\!\! v(y+z,\omega)
    =
\int_0^1 \langle \Phi(\tau_{tx+(1-t)y+z}\,\omega),x-y\rangle\,\dd
t\\
&=& \int_0^1 \langle
\Phi(\tau_{tx+(1-t)y}(\tau_z\,\omega)),x-y\rangle\,\dd t
    =
v(x,\tau_z\omega)-v(y,\tau_z\omega).\nonumber
\end{eqnarray}
for every $\omega$ in a set $\Omega_{x,y}$ of probability 1. Such
a set does not depend only on $z$, as required to prove the claim,
but also on $x$ and $y$. To overcome this difficulty, we set
\[\Omega_z
=\Omega_0\cap\left(\bigcap_{x,y\in\Q^N}\Omega_{x,y}\right).
\]
Clearly $\PP(\Omega_z)=1$, and by the continuity of
$v(\cdot,\omega)$ we derive that, for any $\omega\in\Omega_z$, the
equality \eqref{eq incr staz} now holds for all $x,y\in\R^N$, as
required. This ends the proof.
\end{dimo}

\ \\
\indent We make use of the results presented above to prove
Theorems \ref{ascoli}, \ref{teorema mean 0} and \ref{teorema
constant mean}.\\

\noindent{\bf Proof of Theorem \ref{ascoli}.} The scheme of the
proof is similar to that of Theorem 3.8 in \cite{DavSic08}. We
denote by $\Phi_n$ the functions of $(L^\infty(\Omega))^N$
associated to $v_n$ through Proposition \ref{dif}. Since the
$\Phi_n$ are bounded in $(L^2(\Omega))^N$, a sequence made up by
finite convex combinations of them, say $\Psi_k= \sum_{n \geq n_k}
\lambda_n^k \Phi_n$,  converges to some $\Phi$ in
$(L^2(\Omega))^N$. By Lemma \ref{lemma continuita},   up to
extraction of a subsequence
\begin{equation}\label{ascoli2}
    \int_\gamma \Psi_k(\omega) \rightarrow \int_\gamma \Phi(\omega)
\quad\text{a.s. in $\omega$,}
\end{equation}
for any  piecewise $C^1$ curve $\gamma$.  This implies, in
particular,  $\int_\gamma \Phi(\omega) =0$ a.s. in $\omega$ for
any closed curve $\gamma$.  We can thus associate to $\Phi$ a
Lipschitz random function $v$ with stationary increments using
Proposition \ref{teo null rot}, and we can further assume
$v(0,\omega)=0$ a.s. in $\omega$. Let $w_k$ be a sequence of
Lipschitz random functions with stationary increments defined as
in the statement and  set $g_k(\omega)=-w_k(0,\omega)$. Given  a
point $x$ and a piecewise $C^1$  curve $\gamma$ connecting $0$ to
$x$, we have, a.s. in $\omega$,
\begin{eqnarray*}
  w_k(x,\omega) + g_k(\omega) &=& \int_\gamma \Psi_k(\omega) \quad\text{for any $k\in\N$} \\
  v(x,\omega) &=& \int_\gamma \Phi(\omega).
\end{eqnarray*}
Therefore $w_k(\cdot,\omega) + g_k(\omega)$ converges pointwise to
$v(\cdot,\omega)$ a.s. in $\omega$, by \abra{ascoli2}. By
construction, the sequence $w_k(\cdot,\omega) + g_k(\omega)$ is
almost surely equiLipschitz--continuous and locally equibounded.
By Ascoli Theorem, it must indeed locally uniformly converge to
$v(\cdot,\omega)$ a.s. in $\omega$, as it was to be proved.\qed

\ \\

\noindent{\bf Proof of Theorem \ref{teorema mean 0}.} Let $\kappa$
be a Lipschitz constant for $v(\cdot,\omega)$ for every $\omega$.
We start by proving the sublinearity property assuming the
gradient to have vanishing mean. The first step is to show the
existence of a set $\Omega'$ of probability $1$ such that
\begin{equation}\label{mean1bis}
    \lim_{t\to
    +\infty}\frac{v(t x ,\omega)}{t}=0\qquad\hbox{for any $x\in\R^N$ and  $\omega \in \Omega'$.}
\end{equation}
 We fix $x \in \R^N$ and denote by $\Phi \in (L^\infty(\Omega))^N$ the function associated
 with
 $v$ through Proposition \ref{dif}. Then
\[
\frac{v(t x,\omega)-v(0,\omega)}{t}=\ \int_0^t\taglio \langle
\Phi(\tau_{sx}\,\omega), x\rangle
 \,\dd s \qquad\hbox{for every $t>0$,}
\]
a.s. in $\omega$.  By applying Birkhoff Ergodic Theorem to the
function $\omega \mapsto \langle \Phi(\omega), x\rangle$ and the
dynamical system $(\tau_{sx})_{s \in \R}$, we get  that
\[
\lim_{t\to +\infty}\frac{v(t x,\omega)-v(0,\omega)}{t}
\]
does exist a.s. in $\omega$ and is almost surely equal to some
measurable function $k(\omega)$. It is easy to see that, for every
$z\in\R^N$, $k(\tau_z\omega)=k(\omega)$ a.s. in $\omega$. Indeed,
\[
\frac{v(t x,\tau_z\omega)-v(0,\tau_z\omega)}{t}=\frac{v(t
x+z,\omega)-v(z,\omega)}{t}\qquad\hbox{a.s. in $\omega$}
\]
for $v$ has stationary increments, and
\[
\left|\frac{v(t x+z,\omega)-v(z,\omega)}{t}-\frac{v(t
x,\omega)-v(0,\omega)}{t}\right|\leq \frac{2\,\kappa |z|}{t}.
\]
By ergodicity, we derive that $k(\cdot)$ is a.s. constant, say
equal to some $k\in\R$. Using the Dominated Convergence Theorem we
infer
\[
    k
    =
    \lim_{t\to +\infty}\int_\Omega \frac{v(t x,\omega)-v(0,\omega)}{t}\,\dd \PP (\omega)
    =
    \lim_{t\to +\infty}\int_\Omega \int_0^t\taglio \langle
\Phi(\tau_{sx}\,\omega), x\rangle \,\dd s\,\dd \PP (\omega).
\]
By exploiting the fact that $\PP$ is invariant with respect to the
translations $\tau_y$, we get
\[
    \int_\Omega \int_0^t\taglio \langle \Phi(\tau_{sx}\,\omega),
x\rangle \,\dd s\,\dd \PP (\omega)
    =
     \int_0^t\taglio \int_\Omega
\langle \Phi(\tau_{sx}\,\omega), x\rangle \,\dd \PP (\omega) \,\dd
s= \langle \EE (\Phi), x \rangle,
\]
and  the limit relation in \abra{mean1bis} follows for $\EE
(\Phi)=0$ by hypothesis, at least for some set $\Omega_x$ of
probability $1$ depending on $x$. We then exploit the Lipschitz
character of $v$ to see that \abra{mean1bis} holds with $\Omega' =
\cap_{x\in\Q^N} \Omega_{x}$. We pick $\omega \in \Omega'$; the
family of functions $y \mapsto \frac {v(t y,\omega)}t$, $t \in
\R_+$,  are equibounded and equiLipschitz  continuous, for $y$
varying in $\partial B_1$, and so it uniformly converges  to $0$,
as $t \rightarrow + \infty$, by Ascoli Theorem and
\abra{mean1bis}. Accordingly, given $\eps >0$, we find
\[ \frac{|u(x)|}{|x|} = \frac{|u( |x| \frac x{|x|})|}{|x|} < \eps\]
for $|x|$  large enough, as claimed. \smallskip\\
\indent We proceed to prove the converse implication. By Birkhoff
Theorem we can find an $\omega$  for which the convergence
\abra{mean1} takes place and
\[
\EE( \Phi_1) = \lim_{R\to +\infty}\quad \int_{[-R,R]^N}{\kern
-14.5 mm -} \quad \partial_{x_1}
    v(x,\omega)\,\dd x. \]
Let us denote a point $x$ in $\R^N$ by
$(x_1,y)\in\R\times\R^{N-1}$, we have
\begin{eqnarray*}
&&\left|\int_{[-R,R]^N}{\kern -14.5 mm -} \ \partial_{x_1}
    v(x,\omega)\,\dd x\right|
    =
    \left|\int_{[-R,R]^{N-1}}{\kern -18 mm -}
    \qquad\quad\ \ \left( \int_{-R}^{R}{\kern -7 mm -} \quad\partial_{x_1}
    v(x_1,y,\omega)\,\dd x_1 \right)\dd y\right|\\
&&\qquad\qquad=
    \left|\int_{[-R,R]^{N-1}}{\kern -17.5 mm -}
    \qquad\quad\ \ \left(\frac{v(R,y,\omega)-v(-R,y,\omega)}{R}\right)
    \,\dd y\right|
    \leq
    2\,\max_{[-R,R]^N}\frac{|v(x,\omega)|}{R},
\end{eqnarray*}
which implies $\EE(\Phi_1)=0$ in force of  the assumption.
Similarly we show that $\EE(\Phi_i)$ vanishes for any $i= 2,
\cdots N$.\qed

\ \\

\noindent{\bf Proof of Theorem \ref{teorema constant mean}.} Let
us fix $x$, $y$ in $\R^N$ and denote by $\gamma$ the segment $t y
+(1- t) x$, $t \in [0,1]$, and by $Q$ the vector $\EE(\Phi)$,
where $\Phi$ is the function of $(L^\infty(\Omega))^N$ associated
with $v$ through Proposition \ref{dif}. Using Robin's Theorem, as
in Proposition \ref{dif}, we get
\[
v(y,\omega)-v(x,\omega)=\int_0^1 \langle \Phi(\tau_{\gamma(t)}
\omega),y -x \rangle\,\dd t \quad\text{a.s. in $\omega$,}
\]
and by integrating on $\Omega$
\begin{eqnarray}\label{eq rel}
\EE( v(y,\cdot)- v(x,\cdot))
    =
    \int_\Omega\left(\int_0^1 \langle \Phi(\tau_{\gamma(t)}
\omega),y -x \rangle\,\dd
    t\right)\,\dd \PP\qquad\quad\\
    =
\int_0^1 \langle \int_\Omega \Phi(\tau_{\gamma(t)} \omega) \, \dd
\PP,y -x \rangle\ \dd t
    =
    \langle Q,y-x\rangle.\qquad\ \nonumber
\end{eqnarray}
Now, if $v$ has gradient with vanishing mean , i.e. if $Q=0$, then
\abra{ciao} follows,  conversely, if \abra{ciao}  holds, then it
is enough to take $y-x=Q$ in \eqref{eq rel} to get $Q=0$.\qed
\medskip

\medskip
\end{section}
\end{appendix}


\begin{thebibliography}{99}

\small{
%\bibitem{AmBeVe98} \textsc{M. Amar, G. Bellettini, S. Venturini,}
%Integral representation of functionals defined on curves of
%$W^{1,p}$. {\em Proc. Roy. Soc. Edinburgh,} {\bf 128A} (1998),
%193--217.

%\bibitem{AmAsBu89} \textsc{L. Ambrosio, O. Ascenzi, G. Buttazzo,} {Lipschitz Regularity
%for Minimizers of Integral Functionals with Higly discontinuous
%Integrands}. {\em J. Math. Anal. Appl.} {\bf 142}, no. 2 (1989),
%301--316.

\bibitem{BCD97} \textsc{M. Bardi, I. Capuzzo Dolcetta,}
{Optimal control and viscosity solutions of
Hamilton--Jacobi--Bellman equations.} With appendices by Maurizio
Falcone and Pierpaolo Soravia.  Systems \& Control: Foundations \&
Applications. Birkh\"auser Boston, Inc., Boston, MA, 1997.

\bibitem{Ba94} \textsc{G. Barles,}
{Solutions de viscosit\`e  des \'equations de  Hamilton--Jacobi}.
Math\'{e}matiques \& Applications, { 17}.  Springer--Verlag,
Paris, 1994.

\bibitem{Bu} \textsc{D. Burago, Y. Burago, S. Ivanov,} {A course in
Metric Geometry}. Graduate Studies in Mathematics. {33}. AMS,
Providence, 2001.

%
%\bibitem{BR} \textsc{G. Barles, J.M.Roquejoffre,} Large time behaviour of fronts governed by
% eikonal equations.  {\em Interfaces Free Bound.} {\bf 5} (2003),  83--102.

%\bibitem{BrDav05} \textsc{A. Briani, A. Davini,}
%{Monge solutions for discontinuous Hamiltonians}.  {\em ESAIM
%Control Optim. Calc. Var.}  {\bf 11},  no. 2  (2005), 229--251.

\bibitem{BuGiHi98} \textsc{G. Buttazzo, M. Giaquinta, S. Hildebrandt,} One--dimensional
variational problems. An introduction. Oxford Lecture Series in
Mathematics and its Applications, {\bf 15}. The Clarendon Press,
Oxford University Press, New York, 1998.

%
%\bibitem{CCKS} \textsc{L. Caffarelli, M. G. Crandall, M.  Kocan, A.
%Swiech,} {On viscosity solutions of fully nonlinear equations with
%measurable ingredients.} {\em Comm. Pure Appl. Math.} {\bf 49},
%no. 4 (1996), 365-397.
%
%\bibitem{CaSic04} \textsc{F. Camilli, A. Siconolfi,} Time--dependent
%measurable Hamilton--Jacobi equations. {\em  Comm. Partial
%Differential
%  Equations}  {\bf 30 },  no. 4-6   (2005), 813--847.

\bibitem{CS2} \textsc{F. Camilli, A. Siconolfi,}
Effective Hamiltonian and homogenization of  measurable  Eikonal
equations. {\em  Arch. Ration. Mech. Anal.}  {\bf 183}  (2007),  no. 1,
1--20.

\bibitem{CaVa77} \textsc{G. Castaing, M. Valadier,} Convex Analysis and Measurable
Multifunctions. Springer--Verlag, Berlin, 1977.


%\bibitem{cesari} \textsc{L. Cesari}, Optimization Theory and Applications.
%Springer--Verlag, New York, 1983.

\bibitem{Cl} \textsc{F.H. Clarke,} Optimization and nonsmooth analysis. Wiley, New York,
1983.

%\bibitem{ClVi85} \textsc{F.H. Clarke, R.B. Vinter,} Regularity properties of solutions to the
%basic problem in the Calculus of Variations. {\em Trans. Amer.
%Math. Soc.} {\bf 289} (1985), 73--98.


\bibitem{CoIt00} \textsc{G. Contreras, R. Iturriaga,} Global Minimizers of
Autonomous Lagrangians. 22nd Brazilian Mathematics Colloquium,
IMPA, Rio de Janeiro, 1999.

\bibitem{D1-06} \textsc{A. Davini,} Bolza Problems with
discontinuous Lagrangians and Lipschitz continuity of the value
function.  {\em SIAM J. Control Optim.}  {\bf 46}  (2007),  no. 5,
1897--1921.

\bibitem{DavSic05} \textsc{A. Davini, A. Siconolfi,} A generalized dynamical approach
to the large time behavior of solutions of Hamilton--Jacobi
equations. {\em SIAM J. Math. Anal.}, Vol. 38, no. 2 (2006),
478--502.

\bibitem{DavSic08}
\textsc{A. Davini, A. Siconolfi,} Exact and approximate correctors
for stochastic Hamiltonians: the $1$--dimensional case,  {\em
Preprint} (2008) (available at {\tt
http://cvgmt.sns.it/cgi/get.cgi/papers/davsic08/}).

\bibitem{DS1-08}
\textsc{A. Davini, A. Siconolfi,} On a random analog of
Aubry--Mather Theory for stationary ergodic Hamiltonians. Preprint
(2008).

\bibitem{Ev89} \textsc{L.C. Evans,} The perturbed test function method for
viscosity solutions of nonlinear PDE. {\em Proc. Roy. Soc.
Edinburgh Sect. A}  {\bf 111}  (1989),  no. 3-4, 359--375.

\bibitem{Ev92} \textsc{L.C. Evans,} Periodic homogenisation of certain fully
nonlinear partial differential equations.  {\em Proc. Roy. Soc.
Edinburgh Sect. A}  {\bf 120}  (1992),  no. 3-4, 245--265.


\bibitem{FatFig}
\textsc{A. Fathi, A. Figalli,} Optimal transportation on non
compact manifolds. {\em Israel J. Math.}, to appear.

\bibitem{FatMad}
\textsc{A. Fathi, E. Maderna,} Weak KAM theorem on non compact
manifolds.  {\em NoDEA Nonlinear Differential Equations Appl.}
{\bf 14} (2007),  no. 1-2, 1--27.

\bibitem{FaSic03} \textsc{A. Fathi, A. Siconolfi}, PDE aspects of
Aubry--Mather theory for continuous convex Hamiltonians. {\em
Calc. Var. Partial Differential Equations} {\bf 22},  no. 2
(2005) 185--228.

%\bibitem{ioffe} \textsc{A.D. Ioffe,} On the lower semicontinuity of integral functionals
%I. {\em SIAM J. Control Optim.} {\bf 15} (1977), 521--538.

\bibitem{JiKoOl} \textsc{V.V. Jikov, S.M. Kozlov, O.A.
Oleinik}, Homogenization of differential operators and integral
functionals. Translated from the Russian by G.A. Yosifian.
Springer-Verlag, Berlin, 1994.


\bibitem{Li82} \textsc{P. L. Lions}, Generalized  solutions of Hamilton
Jacobi equations. Research Notes in Mathematics, 69. Pitman
(Advanced Publishing Program), Boston, Mass.-London, 1982.

\bibitem{LiSou03} \textsc{P.L. Lions, P.E. Souganidis,}
Correctors for the homogenization of Hamilton-Jacobi equations in
the stationary ergodic setting. {\em Comm. Pure Appl. Math.} {\bf
56} (2003),  no. 10, 1501--1524.

\bibitem{Molchanov} \textsc{I. Molchanov}, Theory of random sets. Probability and its
Applications (New York). Springer-Verlag London, Ltd., London,
2005.

%\bibitem{NR} \textsc{G. Namah, J.M.Roquejoffre,} The "hump" effect in solid
%propellant combustion. {\em Interfaces Free Bound.} {\bf 2}
%(2000), 449--467.
%
%\bibitem{olech} \textsc{C. Olech,} Weak lower semicontinuity of integral functionals.
%Existence theorem issue. {\em J. Optim. Theory Appl.} {\bf 19}
%(1976), 3--16.
%
%\bibitem{O} \textsc{D. Ostrov,} Solutions of Hamilton-Jacobi equations and
%scalar conservation laws with discontinuous space-time dependance.
%{\em J. Diff. Eq.} {\bf 182} (2002), 51--77.

\bibitem{ReTa00} \textsc{F. Rezakhanlou, J. E. Tarver,} Homogenization for
stochastic Hamilton-Jacobi equations.  {\em Arch. Ration. Mech.
Anal.} {\bf 151}  (2000),  no. 4, 277--309.

%\bibitem{R70} ?? \textsc{R.T.  Rockafellar,} Convex Analysis.
%Princeton Mathematical Series, No. 28. Princeton University Press,
%1970.
%
%\bibitem{RoWe}??  \textsc{R.T.  Rockafellar, R. J.-B. Wets,} Variational Analysis.
%Grundlehren der Mathematischen Wissenschaften [Fundamental
%Principles of Mathematical Sciences], 317. Springer-Verlag,
%Berlin, 1998.

%\bibitem{So03} \textsc{P. Soravia,} {Boundary value problems for
%Hamilton-Jacobi equations with discontinuous Lagrangian}. {\em
%Indiana Univ. Math. J.}, {\bf 51}, no. 2 (2002),  451--477.

\bibitem{Souga99} \textsc{P.E. Souganidis},
Stochastic homogenization of Hamilton-Jacobi equations and some
applications.  {\em Asymptot. Anal.}  {\bf 20}  (1999),  no. 1,
1--11.

\bibitem{Souga-personal} \textsc{P.E. Souganidis}, personal
communication.

%\bibitem{St03} \textsc{T. Str\"omberg,}
%On viscosity solutions of irregular Hamilton-Jacobi equations.
%{\em Arch. Math.} (Basel) {\bf 81},  no. 6 (2003), 678--688.
%
%\bibitem{To15} \textsc{L. Tonelli,} Sur une m\'eth\`ode directe du calcul de variations.
%{\em Rend. Circ. Mat. Palermo} {\bf 39} (1915), 223--264.
%
%\bibitem{To21} \textsc{L. Tonelli,} Fondamenti di Calcolo delle Variazioni. Vol. 1 (1921), Vol. 2
%(1923), Zanichelli, Bologna.

 }


\end{thebibliography}
\end{document}